\documentclass[12pt]{book}
\usepackage{amsmath}
\usepackage{amssymb}
\usepackage{amsfonts}
\usepackage{amsthm}
\usepackage{mathrsfs}
\usepackage{amssymb,amsmath,amscd,graphicx,latexsym}
\usepackage{oldgerm}
\usepackage{enumerate}
\usepackage[latin1]{inputenc}
\usepackage[T1]{fontenc}
\usepackage{epic,eepic,epsfig}
\usepackage{float}
\usepackage{fancyhdr}
\usepackage{graphicx}
\usepackage{amstext}

\renewcommand{\a}{\alpha}
\renewcommand{\b}{\beta}
\renewcommand{\d}{\delta}
\newcommand{\D}{\Delta}

\newcommand{\f}{\frac}
\newcommand{\g}{\gamma}
\newcommand{\G}{\Gamma}
\renewcommand{\l}{\lambda}
\renewcommand{\L}{\Lambda}
\newcommand{\bea}{\begin{eqnarray}}
\newcommand{\eea}{\end{eqnarray}}
\newcommand{\bna}{\begin{eqnarray*}}
\newcommand{\ena}{\end{eqnarray*}}

\renewcommand{\O}{\Omega}
\renewcommand{\le}{\left}
\newcommand{\ri}{\right}
\newcommand{\s}{\sigma}
\renewcommand{\th}{\theta}
\newcommand{\ve}{\varepsilon}
\newcommand{\vp}{\varphi}

\def\ligne{\noindent \rule{\textwidth}{0.5mm} \\}

\def\leq{\leqslant}
\def\geq{\geqslant}

\begin{document}
\openup 1 \jot

\thispagestyle{empty}
\begin{minipage}{25mm}
\vspace{1mm}
\includegraphics[width=35mm]{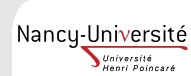}
\end{minipage}
\hfill 
\parbox[c]{10cm}
{\setlength{\arraycolsep}{0cm}\flushright{$\begin{array}{r}
     \mbox{UFR S.T.M.I.A.} \\
     \mbox{\'Ecole Doctorale IAE $+$ M} \\
     \mbox{Universit\'e Henri Poincar\'e - Nancy I} \\
     \mbox{D.F.D. Math\'ematiques}
    \end{array}$}}
\ligne
\vspace{-0.6cm}
\begin{center}
\bf{TH\`{E}SE pr\'esent\'ee par}
\end{center}

\begin{center}
\bf{Yan QU}
\end{center}

\begin{center}
 \bf {pour l'obtention du}
\end{center}
\begin{center} 
\bf{Doctorat de l'Universit\'e Henri Poincar\'e}
\\ Sp\'ecialit\'e : Math\'ematiques. 
\end{center}

\ligne

\begin{center}
\Large{\bf Probl\`{e}mes de type Linnik

pour les fonctions L de formes automorphes.}
\end{center}
\vspace{0.2cm}

\ligne
\vspace{-0.8cm}
\begin{center}
Soutenue publiquement le 2 D\'ecembre 2008
\end{center}

{\bf{Composition du jury :}}

\small{
\begin{tabbing}
\kill\=\hspace{44mm}\=\hspace{40mm}\=\hspace{7cm}\=\kill
\> Emmanuel KOWALSKI     \> Rapporteur  \> Professeur, ETH Zurich\\
\> Yuk-Kam LAU    \> Examinateur  \> Assistant Professeur, Universite Hong Kong\\
\> Jianya LIU             \> Co-directeur de These  \> Professeur, Universite Shandong\\
\> Gerald TENENBAUM            \> President  \> Professeur, Universite Henri Poincare\\
\> Jie WU            \> Co-directeur de These  \> CR au CNRS, HDR, Universite Henri Poincare\\
\end{tabbing}
}

\ligne
\vspace{-1cm}
\begin{center}
\small{Institut \'Elie Cartan Nancy CNRS UMR 9973\\
Facult\'e des Sciences - B.P.239 - 54506 Vandoeuvre-l\`{e}s-Nancy CEDEX}
\end{center}








\tableofcontents

\chapter{Introduction and statement of results}
\section{Linnik's original problem}

In view of Dirichlet's theorem that there are infinitely many primes in the arithmetic
progression
$n\equiv l\,(\bmod \,q)$ with $(q,l)=1$, it is a natural question how big the least prime is, 
denoted by $P(q,l)$, in this arithmetic progression.
Linnik \cite{Lin1} \cite{Lin2} proved that there is an absolute
constant $\ell>0$ such that
$$
P(q,l)\ll q^{\ell},
$$
and this constant $\ell$ was named after him. Since then, a number of authors
have established numerical values for Linnik's
constant $\ell$, while the best result known is
$\ell=5.5$
by Heath-Brown \cite{HeaBro}.
We remark that these results depend on, among other things, numerical
estimates concerning
zero-free regions and the Deuring-Heilbronn phenomenon of Dirichlet $L$-functions.
Under the Generalized Riemann Hypothesis (GRH in brief)
for Dirichlet $L$-functions, the above bounds can be improved to
\bea\label{LinGRH/1}
P(q,l)\ll \vp^2(q)(\log q)^2.
\eea
The conjectured bound is
\bea\label{LinCONJ/1}
P(q,l)\ll_\ve q^{1+\ve}
\eea
for arbitrary $\ve>0$, and this is a consequence of GRH
and another conjecture concerning the universality of the distribution of
nontrivial zeros for Dirichlet $L$-functions.
The conjectured bound (\ref{LinCONJ/1}) is the best possible save the $\ve$
in the exponent. In fact, a trivial lower bound for $P(q,l)$ is
\bea\label{LinTRI/1}
\max_l P(q,l)\geq \{1+o(1)\}\vp(q)\log q.
\eea
Linnik's problem is a rich resource for further mathematical
thoughts, and there are a number of problems that can be formulated
in a similar manner.

\section{A Linnik-type problem for classical modular forms}

\noindent
Let $f$ be a normalized Hecke
eigenform that is a new form of level $N$ of even integral weight $k$ on $\G_0(N)$.
Recall that normalized means that the first Fourier coefficient $\l_f(1)=1$,
and new form means that $N$ is the exact level of $f$, and in this case the Fourier coefficients are equal to
the Hecke eigenvalues.
It also follows that, for this $f$, its Fourier coefficients
$
\{\l_f(n)\}_{n=1}^\infty
$
are real. Applying a classical theorem of Landau, one shows that
the sequence $\{\l_f(n)\}_{n=1}^\infty$ must have
infinitely many sign changes, i.e. there are
infinitely many $n$ such that $\l_f(n)>0$, and there are infinitely
many $n$ such that $\l_f(n)<0$. In view of this result, a reasonable question
to ask is:

{\it Is it possible to obtain a bound on the first sign change, say,
in terms of $k$ and $N$? }

This question is similar to Linnik's problem in nature, and it is named as Linnik-type
in this thesis. In general, this seems to be a difficult question.

In a very special case, this has been considered in
Siegel \cite{Sie}; but, in general, developments have been
achieved only quite recently. In the case $N=1$, sign changes of
the $\l_f(p)$ where $p$ goes over primes have been considered by Ram Murty \cite{Mur}.
Kohnen and Sengupta \cite{KohSen} have shown that the first sign change of $\l_f(n)$
happens for some $n$ with
\bea\label{KohSen/1}
n\ll kN\exp\le(c\sqrt{\f{\log N}{\log\log(3N)}}\ri)(\log k)^{27}, \quad (n,N)=1,
\eea
where $c>2$ is a constant and the $\ll$-constant is absolute. Note that it is
natural to assume that $(n,N)=1$, since the eigenvalues $\l_f(p)$ with $p|N$
are explicitly known by the Atkin-Lehner theory.
Recently, Iwaniec, Kohnen, and Sengupta \cite{IwaKohSen} proved that there is
some $n$ with
\bea\label{IwaKohSen/1}
n\ll (k^2N)^{29/60}, \quad (n,N)=1,
\eea
such that $\l_f(n)<0$.

\medskip

This result is sharp indeed; to see this, let us point out that
the convexity bound
\bea\label{convL12f/1}
L(1/2+it,f)\ll (k^2N)^{1/4+\ve}
\eea
of the automorphic $L$-function $L(1/2+it,f)$ gives, instead of (\ref{IwaKohSen/1}),
the weaker bound
\bea\label{consCONV/1}
n\ll (k^2N)^{1/2+\ve}, \quad (n,N)=1.
\eea
The uniform subconvexity bound
\bea\label{simSCB/1}
L(1/2+it,f)\ll (k^2N)^{29/120}
\eea
would prove (\ref{IwaKohSen/1}), but no result of this quality is known.
The best known uniform subconvexity bound like (\ref{simSCB/1}) is
due to Michel and Venkatesh \cite{MicVen}, which states that
\bea\label{simSCB/MV}
L(1/2+it,f)\ll (k^2N)^{1/4-\d},
\eea
where $\d$ is some positive constant not specified.
Iwaniec, Kohnen, and Sengupta \cite{IwaKohSen}
manage to establish (\ref{IwaKohSen/1}) without appealing to (\ref{simSCB/1});
instead, they use the arithmetic properties of $\l_f(n)$, the Ramanujan conjecture
proved by Deligne, and sieve methods.

A more precise question to ask is: how long is the sequence of Hecke eigenvalues that keep the same sign. To
measure the length of the sequences, define
\bea\label{def/C+-/1}
{\mathscr N}_f^+(x)=\sum_{\substack{n\leq x, \, (n,N)=1 \\
\l_f(n)>0}} 1, \eea and define ${\mathscr N}_f^-(x)$ similarly by replacing the condition $\l_f(n)>0$ under the summation by
$\l_f(n)<0$. Kohnen, Lau, and Shparlinski \cite{KohLauShp} prove that, if $f$ is a new form, then
\bea\label{KLS/17/1}
{\mathscr N}_f^{\pm}(x)\gg_f \f{x}{\log^{17} x},
\eea
where the implied constant depends on the form
$f$. Recently, Wu \cite{Wu} reduces the $17$ in the logarithmic exponent to $1-1/\sqrt{3}$, as an
simple application of his estimates on power sums of Hecke eigenvalues. Still more recently, Lau and Wu
\cite{LauWu} manage to completely get rid of the logarithmic factor in (\ref{KLS/17/1}), getting
\bea\label{KLS/0/1}
{\mathscr N}_f^{\pm}(x)\gg_f x,
\eea
where the implied constant depends on the form $f$. Obviously,
this is the best possible result concerning the order of magnitude of $x$.

These materials form our Chapter 2.

\section{A Linnik-type problem for Maass forms}

\noindent
In Chapter 3, we go on to study Linnik-type problem for Maass forms.
Let $f$ be a normalized Maass
eigenform that is a new form of level $N$ on $\G_0(N)$.
Then, similarly, its Fourier coefficients
$
\{\l_f(n)\}_{n=1}^\infty
$
are real. Applying Landau's theorem as in the holomorphic case,
one shows that the sequence $\{\l_f(n)\}_{n=1}^\infty$ must have
infinitely many sign changes. Therefore, one may formulate a
Linnik-type problem for this normalized Maass
eigenform $f$.

In this direction, we prove the following theorem.

\noindent
{\bf Theorem 3.10.}
{\it Let $f$ be a normalized Maass new form of
level $N$ and Laplace eigenvalue $1/4+\nu^2$.
Then there is some $n$ satisfying
\bea\label{n<</MAASS/1}
n\ll ((3+|\nu|)^2N)^{1/2-\delta}, \quad (n,N)=1,
\eea
such that $\l_f(n)<0$, where $\delta$ is a positive absolute constant. }

This is proved by using, among other things, the uniform
subconvexity bound (\ref{simSCB/MV})
of Michel and Venkatesh \cite{MicVen}, and
this explains why we are not able to get an acceptable numerical value for $\delta$.
However, the bound (\ref{simSCB/MV}) alone is not enough to establish
(\ref{n<</MAASS/1}); some combinatorial and analytic
arguments are also necessary for (\ref{n<</MAASS/1}).
The method in Iwaniec, Kohnen, and Sengupta \cite{IwaKohSen} does not
work here; one of the reasons is that for Maass forms
the Ramanujan conjecture is still open, and hence the sieves in \cite{IwaKohSen}
do not apply.

\section{A Linnik-type problem for automorphic $L$-functions}

The Linnik-type problem considered before can be further generalized to that
for automorphic $L$-functions. This is done in Chapter 4.

To each irreducible unitary cuspidal representation
$\pi=\otimes \pi_p$ of $GL_m({\mathbb A}_{\mathbb Q})$, one can attach a global $L$-function
$L(s,\pi)$, as in Godement and Jacquet \cite{GodJac}, and Jacquet and
Shalika \cite{JacSha1}. For $\s=\Re s>1$, $L(s,\pi)$ is defined by
products of local
factors
\bea\label{Lspi/1}
L(s,\pi)
=\prod_{p<\infty} L_p (s, \pi_p),
\eea
where
\bea\label{Lpspip/1}
L_p(s, \pi_p )=\prod_{j=1}^m
\le(1-\f{\alpha_\pi(p,j)}{p^{s}}\ri)^{-1};
\eea
the complete $L$-function $\Phi(s,\pi)$ is defined by
\bea\label{PHIspi/1}
\Phi(s,\pi)=L_\infty(s,\pi_\infty)L(s,\pi),
\eea
where
\bea\label{L8spi8/1}
L_\infty(s, \pi_\infty)=
\prod_{j=1}^m
\Gamma_{\mathbb R}(s+\mu_{\pi}(j))
\eea
is the Archimedean local factor. Here
\bea\label{GammaR}
\Gamma_{\mathbb R}(s) =\pi^{-s/2}\Gamma \le(\f{s}{2}\ri),
\eea
and $\{\alpha_\pi(p,j)\}_{j=1}^m$ and $\{\mu_{\pi} (j)\}_{j=1}^m$
are complex numbers associated with $\pi_p$ and $\pi_\infty$,
respectively, according to the Langlands correspondence.
The case $m=1$ is classical; for $m\geq 2$, $\Phi(s,\pi)$ is
entire and satisfies a functional equation.

It is known from Jacquet and Shalika \cite{JacSha1}
that the Euler product for $L(s,\pi)$ in
(\ref{Lspi/1}) converges absolutely for $\s>1$.
Thus, in the half-plane $\s>1$, we may write
\bea\label{Lspi=Dir/1}
L(s,\pi)=\sum_{n=1}^\infty \f{\l_\pi(n)}{n^s},
\eea
where
\bea\label{lpin=/1}
\l_\pi(n)=\prod_{p^\nu\| n}
\bigg\{\sum_{\nu_1+\cdots +\nu_m=\nu} \a_\pi(p,1)^{\nu_1}\cdots \a_\pi(p,m)^{\nu_m} \bigg\}.
\eea
In particular,
\bea\label{lpin=PAR/1}
\l_\pi(1)=1,
\qquad
\l_\pi(p)=\a_\pi(p,1)+\cdots +\a_\pi(p,m).
\eea
It also follows from work of Shahidi \cite{Sha1}, \cite{Sha2}, \cite{Sha3}, and \cite{Sha4}
that the complete $L$-function
$\Phi(s,\pi)$ has an analytic continuation to
the whole complex plane and satisfies the functional equation
$$
\Phi(s,\pi) =\ve(s,\pi)
\Phi(1-s,\tilde{\pi})
$$
with
$$
\ve(s,\pi)=
\ve_\pi N_{\pi}^{1/2-s},
$$
where $N_{\pi}\geq 1$ is an integer called the arithmetic conductor of $\pi$, $\ve_\pi$
is the root number satisfying $|\ve_\pi|=1$,
and $\tilde{\pi}$ is the representation contragredient to $\pi$.

By an argument similar to the case of holomorphic forms or Maass forms,
it is possible to establish the following theorem of
infinite sign changes.

\noindent
{\bf Theorem 4.13.} {\it Let  $m\geq 2$ be an inetger and 
let $\pi$ be an irreducible unitary cuspidal representation
for $GL_m(\mathbb A_\mathbb Q)$ such that $\l_\pi(n)$ is real for all $n\geq 1$. 
Then the sequence
$
\{\l_f(n)\}_{n=1}^\infty
$
has infinitely many sign changes, i.e. there are
infinitely many $n$ such that $\l_f(n)>0$, and there are infinitely
many $n$ such that $\l_f(n)<0$.
}

Iwaniec and Sarnak \cite{IwaSar} introduced the analytic
conductor of $\pi$. It is a function over the reals given by
\bea\label{AnalCod/1}
Q_\pi(t)=N_\pi \prod_{j=1}^m (3+|t+\mu_\pi(j)|),
\eea
which puts together all the important parameters for $\pi$. The
quantity
\bea\label{Cod/1}
Q_\pi=Q_\pi(0)=N_\pi \prod_{j=1}^m (3+|\mu_\pi(j)|)
\eea
is also important, and it is named the conductor of $\pi$.

We may therefore formulate a Linnik-type problem for
$\{\l_f(n)\}_{n=1}^\infty$, and
the first sign change is measured by
the conductor $Q_\pi$ of $\pi$. Our result in this direction is as follows.

\noindent
{\bf Theorem 4.15.}
{\it Let  $m\geq 2$ be an inetger and let $\pi$ be an irreducible unitary cuspidal representation
for $GL_m(\mathbb A_\mathbb Q)$.
If $\l_\pi(n)$ is real for all $n\geq 1$, then there is some $n$ satisfying
\bea\label{bound/AUTO/1}
n\ll Q_\pi^{m/2+\ve}
\eea
such that $\l_\pi(n)<0$. The constant implied in (\ref{bound/AUTO/1})
depends only on $m$ and $\ve$. In particular, the result is true
for any self-contragredient representation $\pi$. }

Proof of this theorem is quite different from that
of Theorem 3.10.
One of the most principal difficulties is that there is no relation of Hecke-type
in this current general case of $\pi$ being
irreducible unitary cuspidal representation, as in classic modular forms or Maass form cases
(see (2.5) and (3.11)).
These difficulties
are overcome by, among other things,
new analytic properties of $L(s,\pi)$ due to Harcos \cite{Har},
an inequality due to Brumley \cite{Bru}, as well as important combinatorial properties
of the sequence $\{\l_f(n)\}_{n=1}^\infty$, established in Lemma 4.12.
Clearly this elegant inequality is of independent interest
and we believe that it will find other applications.

\section{Automorphic prime number theorem and a problem of Linnik's type}

To each irreducible unitary cuspidal representation
$\pi=\otimes \pi_p$ of $GL_m({\mathbb A}_{\mathbb Q})$, one can attach a global $L$-function
$L(s,\pi)$ as in \S 1.4.
Then, one can link $L(s,\pi)$ with primes by taking logarithmic differentiation
in (\ref{Lpspip/1}), so that for $\s>1$,
\bea\label{1/DEFLa}
\frac{\rm d}{{\rm d}s} \log L(s,\pi)
=
-\sum_{n=1}^\infty
\frac{\L(n)a_\pi(n)}{n^s},
\eea
where $\L(n)$ is the von Mangoldt function, and
\bea\label{apipk/AUTO/1}
a_\pi(p^k)=
\sum_{j=1}^m
\alpha_\pi(p,j)^k.
\eea
The prime number theorem for $L(s,\pi)$ concerns the asymptotic
behavior of the counting function
\bna
\psi(x,\pi)=\sum_{n\leq x}\L(n)a_\pi(n),
\ena
and a special case of it
asserts that, if $\pi$ is an irreducible unitary cuspidal
representation of $GL_m({\mathbb A}_{\mathbb Q})$ with $m\geq 2$, then
\bea\label{psixpi/AUTO/1}
\psi(x,\pi)\ll \sqrt{Q_\pi} \cdot x \cdot \exp\bigg(-\f{c}{2m^4}\sqrt{\log x}\bigg)
\eea
for some absolute positive constant $c$, where the implied constant is absolute.
In Iwaniec and Kowalski
\cite{IwaKow}, Theorem 5.13, a prime number theorem is proved for general
$L$-functions satisfying necessary axioms, from which (\ref{psixpi/AUTO/1})
follows as a consequence.

In this chapter, we first investigate the influence of GRH
on $\psi(x,\pi)$. It is known that, under GRH, (\ref{psixpi/AUTO/1}) can be improved to
\bea\label{undGRH/AUTO/1}
\psi(x,\pi)\ll x^{1/2}\log^2 (Q_{\pi}x),
\eea
where the implied constant depends at most on $m$.
But better results are desirable.
In this direction, we establish the following results.

\noindent
{\bf Theorem 5.1.} {\it Let  $m\geq 2$ be an inetger and let $\pi$ be an irreducible unitary cuspidal
representation of $GL_m({\mathbb A}_{\mathbb Q})$. Assume GRH for $L(s,\pi)$. Then
we have
\bna
\psi(x,\pi)\ll x^{1/2}\log^2(Q_{\pi}\log x)
\ena
for $x\geq 2$, except on a set $E$ of finite logarithmic measure, i.e.
\bna
\int_{E} \f{{\rm d}x}{x}< \infty.
\ena
The constant implied in the
$\ll$-symbol depends
at most on $m$.}

\noindent
{\bf Theorem 5.2.} {\it Let $m\geq 2$ be an inetger and let $\pi$ be an irreducible unitary cuspidal
representation of $GL_m({\mathbb A}_{\mathbb Q})$. Assume GRH for $L(s,\pi)$. Then
\bna
\int^X_2 |\psi(x,\pi)|^2\f{{\rm d}x}{x}\ll X\log^2Q_{\pi}.
\ena
The constant implied in the
$\ll$-symbol depends
at most on $m$.}

Gallagher \cite{Gal1} was the first to establish a result
like Theorem 5.1, in the classical case $m=1$ for the Riemann
zeta-function. He proved that, under the Riemann Hypothesis for
the classical zeta-function,
$$
\psi(x):=\sum_{n\leq x}\L(n)=x+O\big(x^{1/2}(\log \log x)^2\big)
$$
for $x\geq 2$, except on a set of finite logarithmic measure,
and hence made improvement on the classical estimate error term
$O(x^{1/2}\log^2 x)$ of von Koch \cite{Koc}. In the same paper, Gallagher \cite{Gal1}
also gave short proofs of Cram\'er's conditional estimate (see \cite{Cra1} \cite{Cra2})
\bna
\int_2^X (\psi(x)-x)^2\f{{\rm d}x}{x}\ll X.
\ena
Gallagher's
proofs of the above results make crucial use of his lemma in \cite{Gal2},
which is now named after him.

Our Theorems 5.1-5.2 generalize the above classical results to
the prime counting function $\psi(s,\pi)$ attached to
irreducible unitary cuspidal
representations $\pi$ of $GL_m({\mathbb A}_{\mathbb Q})$ with $m\geq 2$.
Our proofs combine the approach of Gallagher with recent results
of Liu and Ye (\cite{LiuYe1}, \cite{LiuYe3})
on the distribution of zeros of Rankin-Selberg automorphic $L$-functions.

The above Theorem 5.2 states that, under GRH, $|\psi(x,\pi)|$ is of size $x^{1/2}\log Q_\pi$
on average. This can be compared with the next theorem, which gives
the unconditional Omega result that $|\psi(x,\pi)|$ should not be of order lower
than $x^{1/2-\ve}$.

\noindent
{\bf Theorem 5.3.} {\it Let $m\geq 2$ be an integer and let $\pi$ be an irreducible unitary cuspidal
representation of $GL_m({\mathbb A}_{\mathbb Q})$, and $\ve>0$ arbitrary.
Unconditionally,
\bna
\psi(x,\pi)=\Omega(x^{1/2-\ve}),
\ena
where the implied constant depends at most on $m$ and $\varepsilon$.
More precisely,
there exists an increasing sequence $\{x_n\}_{n=1}^\infty$ tending to infinity such that
\bea\label{upperLIM/1}
\lim_{n\to\infty}
\f{|\psi(x_n,\pi)|}{ x_n^{1/2-\ve}}>0.
\eea
}

Note that the sequence $\{x_n\}_{n=1}^\infty$ and the limit in (\ref{upperLIM/1})
may depend on $\pi$.
This result generalizes that for the Riemann zeta-function. It is possible to get
better Omega results like those in Chapter V of Ingham \cite{Ing}.
We remark that, unlike the classical case, in Theorems 5.1-5.3 we do not have
the main term $x$. This is because $L(s,\pi)$
is entire when $m\geq 2$, while $\zeta(s)$ has a simple pole at
$s=1$ with residue $1$.

Connecting with Linnik's
problem for automorphic $L$-functions considered in Chapter 4,
we consider a Linnik-type problem in the sequence
$
\{a_\pi(n)\L(n)\}_{n=1}^\infty,
$
defined as in (\ref{1/DEFLa}) and (\ref{apipk/AUTO/1}). These
are the coefficients in the Dirichlet series expansion for
$-\f{L'}{L}(s,\pi)$ with $\s>1$. As a consequence of Theorem 4.15,
we establish in this direction  the following theorem.

\noindent
{\bf Theorem 5.12.}
{\it Let $m\geq 2$ be an integer and
let $\pi$ be an irreducible unitary cuspidal representation for $GL_m(\mathbb A_\mathbb Q)$.
If all $a_\pi(n)\L(n)$ are real, then $\{a_\pi(n)\L(n)\}_{n=1}^\infty$ changes sign at some $n$ satisfying
\bea\label{bound/AUTO-A/1}
n\ll Q_\pi^{m/2+\ve}.
\eea
The constant implied in (\ref{bound/AUTO-A/1})
depends only on $m$ and $\ve$. In particular, the result is true
for any self-contragredient representation $\pi$. }

These are the materials in Chapter 5.

\section{Selberg's normal density theorem for automorphic $L$-functions}

Under the Riemann Hypothesis for the Riemann zeta-function,
i.e. in the case of $m=1$, Selberg \cite {Sel} proved that
\bea
\int^X_1\{\psi(x+h(x))-\psi(x)-h(x)\}^2 dx=o(h(X)^2 X)
\eea
for any increasing functions $h(x) \leq x$ with
$$\f{h(x)}{\log^2 x}\to \infty,
$$
where as usual,
$$
\psi(x)=\sum_{n\leq x}\L(n).
$$
In Chapter 6, we prove an analogue of this in the case of automorphic $L$-functions.

\medskip
\noindent
{\bf Theorem 6.1.} {\it Let $m\geq 2$ be an integer and let $\pi$ be an irreducible unitary cuspidal
representation of $GL_m({\mathbb A}_{\mathbb Q})$.
Assume GRH for $L(s,\pi)$. We have
\bea\label{Sel/1}
\int^X_1 |\psi(x+h(x),\pi)-\psi(x,\pi)|^2 dx=o(h(X)^2 X),
\eea
for any increasing functions $h(x)\leq x$ satisfying
$$
\f{h(x)}{\log^2 (Q_{\pi}x)}\to \infty.
$$
}

Our Theorem 6.1 generalizes Selberg's result to cases when $m\geq 2$.
It also improves an earlier result of the author \cite{Qu2}
that (\ref{Sel/1})
holds for $h(x)\leq x$ satisfying
$$
\f{h(x)}{x^{\th}\log^2 (Q_\pi x)}\to \infty,
$$
where $\th$ is the bound towards the GRC as explained in Lemma 4.8.
The main new idea is a delicate application of Kowalski-Iwaniec's mean value estimate (cf. Lemma 5.9).
We also need an explicit formula established in Chapter 5 in a more precise form.

Unconditionally, Theorem 6.1 would hold for $h(x)=x^{\b}$
with some constant $0<\b<1$. The exact value of $\b$ depends
on two main ingredients: a satisfactory zero-density estimate for
the $L$-function $L(s,\pi)$, and a zero-free region for $L(s,\pi)$
of Littlewood's or Vinogradov's type.

\chapter{Classical modular forms and a Linnik-type problem}

\noindent
In this chapter, we will review the concept and some basic properties of
classical modular forms. These properties will be used in later chapters of the
thesis. The reader is referred to Iwaniec \cite{Iwa1} for a detailed treatment of
these materials.

\section{Classical modular forms}

\noindent
Let
$$SL_2(\mathbb Z)
:=\le\{
\begin{pmatrix}a & b\\c & d\end{pmatrix} :
a, b, c, d\in {\mathbb Z}, \; ad-bc=\pm 1\ri\}
$$
be the modular group. We restrict our attention
to the Hecke congruence subgroup of level $N$, which is
$$
\Gamma_0(N)
=\le\{\left(
       \begin{array}{cc}
         a & b \\
         c & d \\
       \end{array}
     \right)
\in SL_2(\mathbb Z): N|c\ri\},
$$
where $N$ is a positive integer. In this convention, $\Gamma_0(1)=SL_2(\mathbb Z)$, and
the index of $\Gamma_0(N)$ in the modular group is
$$
\nu(N)=[\G_0(1):\G_0(N)]=N\prod_{p|N}\le(1+\f1p\ri).
$$
The group $\G_0(N)$ acts on the upper half-plane
$$
\mathbb H=\{z:\ z=x+iy, \ y>0\}
$$
by
$$
\g z=\f{az+b}{cz+d}, \quad \g=
\left(
       \begin{array}{cc}
         a & b \\
         c & d \\
       \end{array}
     \right)
\in \G_0(N).
$$
Let $k$ be a positive integer. The space of cusp forms
of weight $k$ and level $N$ is denoted by $S_k(\G_0(N))$;
it is a finite-dimensional Hilbert space with respect to
the inner product
\bna
\langle f,g\rangle=\int_{\G\backslash \mathbb H} f(z)\bar{g}(z)y^k\f{{\rm d} x \, {\rm d} y}{y^2},
\ena
where
$$
{\rm d} \mu:=\f{{\rm d} x \, {\rm d} y}{y^2}
$$
is the invariant measure on $\mathbb H$.

The Hecke operators $\{T_n\}_{n=1}^\infty$ are defined by
\bea\label{Tnfz=}
(T_n f)(z)=\f{1}{\sqrt{n}}\sum_{ad=n}\le(\f{a}{d}\ri)^{k/2}
\sum_{b\,({\rm mod}\,d)}f\le(\f{az+b}{d}\ri).
\eea
It follows that $T_m$ and $T_n$ commute, and for every $n$,
$T_n$ is also self-adjoint on $S_k(\G_0(N))$, i.e.
$$
\langle T_n f, g\rangle=\langle f, T_n g\rangle, \quad (n,N)=1.
$$
Let $\mathcal F=\{f\}$ be an orthonormal basis
of $S_k(\G_0(N))$. We can assume that every $f\in \mathcal F$ is an eigenfunction
for all Hecke operators $T_n$ with $(n,N)=1$; i.e.
there exist complex numbers $\l_f(n)$, such that
\bea\label{Tnf=lfnf}
T_n f=\l_f(n)f, \quad (n,N)=1.
\eea
The eigenvalues $\l_f(n)$ are related to the Fourier
coefficients of $f(z)$ in such a way that the Fourier
series expansion of $f(z)$ now takes the form
\bea\label{fz=fourier}
f(z)=\sum_{n=1}^\infty a_f(n)n^{(k-1)/2}e(nz),
\eea
with
\bea\label{afn=af1}
a_f(n)=a_f(1)\l_f(n), \quad (n,N)=1.
\eea
Note that if $a_f(1)=0$, then all $a_f(n)=0$ for $(n,N)=1$.
Here we have used the stardard notation
$$
e(t) := e^{2\pi it}
\qquad(t\in \mathbb R).
$$

\noindent
{\bf Lemma 2.1.} {\it Let $f$ be an eigenfunction 
for all Hecke operators $T_n$ with $(n,N)=1$, and
$\l_f(n)$ be as in (\ref{Tnf=lfnf}). Then

\vskip -2mm

\mbox{\rm (i)} The Hecke eigenvalues $\{\l_f(n)\}_{n=1}^\infty$ are real;

\vskip -2mm

\mbox{\rm (ii)}  The Hecke eigenvalues are multiplicative
in the following sense:
$$
\l_f(m)\l_f(n)=\sum_{d|(m,n)}\l_f\le(\f{mn}{d^2}\ri),
\quad (n,N)=1.
$$
In particular
\bea\label{lf(p)^2=}
\l_f(p)^2=\l_f(p^2)+1, \quad (p,N)=1.
\eea}

Unfortunately, we cannot deduce from (\ref{afn=af1}) that
$a_f(n)\not=0$ because the condition does not allow us to
control all the coefficients in (\ref{fz=fourier}). However,
for new forms, the following result is true.

\noindent
{\bf Lemma 2.2.} {\it If $f$ is a new form, then
(\ref{Tnf=lfnf}) holds for all $n$. The first coefficient
in the Fourier expansion (\ref{fz=fourier}) does not vanish,
so one can normalize $f$ by setting $a_f(1)=1$. In this case,
$a_f(n)=\l_f(n)$ for all $n$, and hence the Fourier expansion of
$f$ takes the form
\bea\label{fz=fourier/NEW}
f(z)=\sum_{n=1}^\infty \l_f(n)n^{(k-1)/2}e(nz).
\eea
 }

Lemmas 2.1 and 2.2 will be used several times later.

\section{Classical automorphic $L$-functions}

\noindent
We begin with a cusp form which has Fourier expansion as in
(\ref{fz=fourier}) and (\ref{afn=af1}). Define, for $\s>1$,
\bea\label{def/Lsf}
L(s,f)=\sum_{n=1}^\infty \f{\l_f(n)}{n^s}.
\eea
The so-called complete $L$-function is defined as
\bea\label{def/compLsf}
\Phi(s,f)=\pi^{-s}\G\le(\f{s+(k-1)/2}{2}\ri)\G\le(\f{s+(k+1)/2}{2}\ri)L(s,f).
\eea
This complete $L$-function satisfies the functional equation
\bea\label{funeqn}
\Phi(s,f)=\ve_f N^{1/2-s}\Phi(1-s,\bar{f}),
\eea
where $\ve_f$ is a complex
number of modulus $1$. For any new form $f$, we have the following
Euler product for $\Phi(s,f)$.

\noindent
{\bf Lemma 2.3.} {\it If $f$ is a new form, then the functional
equation takes the form
\bea\label{funeqn/NEW}
\Phi(s,f)=\ve_f N^{1/2-s}\Phi(1-s,f), \qquad s\in \mathbb C.
\eea
For $\s>1$,  the function $L(s,f)$ admits the Euler product
\bea\label{EuProd}
L(s,f)=\prod_{p}\le(1-\f{\l_f(p)}{p^s}+\f{\chi_N^0(p)}{p^{2s}}\ri)^{-1},
\eea
where $\chi_N^0$ is the principal character modulo $N$.}

By Lemmas 2.1 and 2.2, all eigenvalues $\l_f(n)$ of a new form $f$ are real.
This explains why on the right-hand side of (\ref{funeqn/NEW}) we write
$\Phi(1-s,f)$ instead of $\Phi(1-s,\bar f)$.
We may further factor the Hecke polynomial in (\ref{EuProd}) into
$$
1-\f{\l_f(p)}{p^s}+\f{\chi_N^0(p)}{p^{2s}}=
\le(1-\f{\a_f(p)}{p^s}\ri)\le(1-\f{\b_f(p)}{p^s}\ri),
$$
where
$$
\begin{cases}
\a_f(p)+\b_f(p)=\l_f(p),
\\\noalign{\vskip 2mm}
\a_f(p)\b_f(p)=\chi^0_N(p).
\end{cases}
$$
Recall that the Ramanujan conjecture asserts that
\bea\label{Ramun/Del}
|\a_f(p)|=|\b_f(p)|=1, \quad (p,N)=1,
\eea
which has been proved by Deligne \cite{Del}. For $\s>1,$ the symmetric square
$L$-function is defined by
\bea\label{defsym2}
L(s,\text{\rm sym}^2f)=L(2s,\chi_N^0)\sum_{n=1}^\infty \f{\l_f(n^2)}{n^s}.
\eea
The Euler product of $L(s,\text{\rm sym}^2f)$ takes the form that, for $\s>1$,
\bea\label{SYMLeuler}
L(s,\text{\rm sym}^2f)
& := &\prod_{p}\le(1-\f{\a_f(p)\a_f(p)}{p^{s}}\ri)^{-1}
\le(1-\f{\a_f(p)\b_f(p)}{p^{s}}\ri)^{-1}
\nonumber\\
&& \times \le(1-\f{\b_f(p)\b_f(p)}{p^{s}}\ri)^{-1}.
\eea

Now suppose $g$ is a new form of level $N'$ and weight $k'$.
The Rankin-Selberg $L$-function of $f$ and $g$ is defined as
\bea\label{defRSL}
L(s,f\otimes g)=L(2s,\chi_N^0\chi_{N'}^0)\sum_{n=1}^\infty \f{\l_f(n)\l_g(n)}{n^s},
\eea
if $[N,N']$ is square-free.
For $\s>1$, the Euler product of $L(s,f\otimes g)$ takes the form
\bea\label{RSLeuler}
L(s,f\otimes g)
&=&\prod_{p}\le(1-\f{\a_f(p)\a_g(p)}{p^{s}}\ri)^{-1}
\le(1-\f{\a_f(p)\b_g(p)}{p^{s}}\ri)^{-1}\nonumber\\
&&\times \le(1-\f{\b_f(p)\a_g(p)}{p^{s}}\ri)^{-1}
\le(1-\f{\b_f(p)\b_g(p)}{p^{s}}\ri)^{-1}.
\eea

These properties are important in our later argument.

\section{Infinite sign changes of Fourier coefficients}

The result in the following seems well known, but we are
not able to give a precise reference where it appeared first.
As a substitute, we refer to the paper \cite{KnoKohPri} for an
extension to quite general subgroups of $SL_2(\mathbb R)$ and a
discussion of related topics.

\noindent
{\bf Proposition 2.4.} \text{(Knopp-Kohnen-Pribitkin \cite{KnoKohPri}).}
{\it Let $f$ be a non-zero
cusp form of even integral weight $k$ on $\G_0(N)$,
and suppose that its Fourier coefficients $a_f(n)$ are real for all
$n\geq 1$. Then  the sequence
$$
\{a_f(n)\}_{n=1}^\infty
$$
has infinitely many sign changes, i.e. there are
infinitely many $n$ such that $a_f(n)>0$, and there are infinitely
many $n$ such that $a_f(n)<0$.
}

It follows from Lemmas 2.1 and 2.2 that, for any new form $f$, all its eigenvalues $\l_f(n)$
are real. We therefore arrive at the following corollary.

\noindent
{\bf Corollary 2.5.} {\it Let $f$ be a normalized Hecke
eigenform that is a new form of level $N$ of even integral weight $k$ on $\G_0(N)$.
Then  the sequence
$$
\{\l_f(n)\}_{n=1}^\infty
$$
has infinitely many sign changes, i.e. there are
infinitely many $n$ such that $\l_f(n)>0$, and there are infinitely
many $n$ such that $\l_f(n)<0$.
}

A more precise question to ask is: how long is the sequence of Hecke eigenvalues
that keep the same sign. To measure the length of the sequences, define
\bea\label{def/C+-}
{\mathscr N}_f^+(x)=\sum_{\substack{n\leq x, \, (n,N)=1\\ \l_f(n)>0}} 1,
\eea
and define ${\mathscr N}_f^-(x)$ similarly by replacing the condition $\l_f(n)>0$
under the summation by $\l_f(n)<0$. Kohnen, Lau, and Shparlinski \cite{KohLauShp}
prove that, if $f$ is a new form, then
\bea\label{KLS/17}
{\mathscr N}_f^{\pm}(x)\gg_f \f{x}{\log^{17} x},
\eea
where the implied constant depends on the form $f$.
Recently, Wu \cite{Wu}
reduces the $17$ in the logarithmic exponent to  $1-1/\sqrt{3}$,
as an simple application of his estimates on power sums of Hecke eigenvalues.
Still more recently, Lau and Wu \cite{LauWu} manage to completely
get rid of the logarithmic factor in (\ref{KLS/17}).

\noindent
{\bf Proposition 2.6.} \text{(Lau-Wu \cite{LauWu}).}
{\it Let $f$ be a normalized Hecke
eigenform that is a new form of level $N$ of even integral weight $k$ on $\G_0(N)$,
and let ${\mathscr N}_f^\pm(x)$ be as in (\ref{def/C+-}). Then
\bea\label{KLS/0}
{\mathscr N}_f^{\pm}(x)\gg_f x,
\eea
where the implied constant depends on the form $f$. }

Obviously, this is the best possible result concerning the order of magnitude of $x$.
The proof applies, among other
things, the $\mathscr B$-free number method. It is also remarked in \cite{LauWu} that their method works well in
other cases, such as forms of half-integral weight.

\section{A Linnik-type problem: the first sign change of Fourier coefficients}

\subsection{Linnik's original problem}
In view of Dirichlet's theorem that there are infinitely many primes in the arithmetic
progression
$n\equiv l(\bmod q)$ with $(q,l)=1$, it is a natural question how big the least prime is, 
denoted by $P(q,l)$, in this arithmetic progression.
Linnik \cite{Lin1} \cite{Lin2} proved that there is an absolute
constant $\ell>0$ such that
$$
P(q,l)\ll q^{\ell},
$$
and this constant $\ell$ was named after him. Since then, a number of authors
have established numerical values for Linnik's
constant $\ell$, while the best result known is
$\ell=5.5$
by Heath-Brown \cite{HeaBro}.
We remark that these results depend on, among other things, numerical
estimates concerning
zero-free regions and the Deuring-Heilbronn phenomenon of Dirichlet $L$-functions.
Under GRH for Dirichlet $L$-functions, the above bounds can be improved to
\bea\label{LinGRH}
P(q,l)\ll \vp^2(q)(\log q)^2.
\eea
The conjectured bound is
\bea\label{LinCONJ}
P(q,l)\ll_\ve q^{1+\ve}
\eea
for arbitrary $\ve>0$, and this is a consequence of GRH
and another conjecture concerning the universality of the distribution of
nontrivial zeros for Dirichlet $L$-functions, as shown in Liu and Ye \cite{LiuYe4}.
The conjectured bound (\ref{LinCONJ}) is the best possible save the $\ve$
in the exponent. In fact, a trivial lower bound for $P(q,l)$ is
\bea\label{LinTRI}
\max_l P(q,l)\geq \{1+o(1)\}\vp(q)\log q.
\eea

The reader is referred to \cite{HeaBro}
for a survey of results concerning Linnik's problem.

\subsection{A Linnik-type problem and a classical result of Siegel}

According to the results in the previous section, a reasonable question
to ask is: {\it Is it possible to obtain a bound on the first sign change, say,
in terms of $k$ and $N$? }
In general, this seems to be a difficult question. For a survey of results in this
direction, see Kohnen \cite{Koh1}.

If $f\not=0$, then, by the valence formula for modular forms,
the orders of zeros of $f$ on the compactified Riemann surface
$$
X_0(N)=\G_0(N)\backslash \mathbb H\cup \mathbb P^{1}(\mathbb Q)
$$
sum up to $\f{k}{12}[\G_0(1):\G_0(N)]$.
Hence there exists a number $n$ in the range
$$
1\leq n\leq \f{k}{12}[\G_0(1):\G_0(N)]
$$
such that $a_f(n)\not=0$. Now if we are optimistic,
then we can expect a sign change in the range
$$
1\leq n\leq \f{k}{12}[\G_0(1):\G_0(N)]+1.
$$

In a very special case, this indeed follows from
work of Siegel \cite{Sie}. To formulate the result,
suppose $k\geq 4$ and denote by $d_k$ the dimension
of the space $M_k(\G_0(1))$ of modular forms of
weight $k$ on $\G_0(1)$. Recall that $d_k$ satisfies the
formula
\bna
d_k=
\left\{
\begin{array}{lll}
[k/12]   & k\equiv 2\,(\bmod 12),
\\\noalign{\vskip 3mm}
[k/12]+1 & \text{otherwise}.
\end{array}
\right.
\ena
Then Siegel showed that, for each $f\in M_k(\G_0(1))$, there are explicitly computable
rational numbers $\{c_n\}_{n=0}^{d_k}$ depending on $k$,
such that
$$
\sum_{n=0}^{d_k}c_na_f(n)=0.
$$
Siegel's explicit expression for $c_n$ implies that,
for $k\equiv 2\,(\bmod\,4)$, all the $c_n$ are strictly
positive. Since a cusp form of weight $k$ on $\G_0(1)$
is determined by its Fourier coefficients $\{a_f(n)\}_{n=0}^{d_k-1}$,
we conclude immediately that, under the assumption
$k\equiv 2\,(\bmod\,4)$, there must be a sign change of $a_f(n)$ in the
range $1\leq n\leq d_k$. Thus, using the formula for $d_k$ above, one sees
that the above optimistic expectation is justified in this special case.

Unfortunately, when $k\equiv 0\,(\bmod\,4)$ or if $N>1$, Siegel's argument
does not work any longer, and therefore we need other ideas.

\subsection{Recent developments and comments}

In this subsection, we will focus on the
case that $f$ is a normalized Hecke eigenform that is a new form of level $N$.
Recall that normalized means that $a_f(n)=1$, and new form means that $N$ is the
exact level of $f$, and in this case the Fourier coefficients are equal to
the Hecke eigenvalues.

In the case $N=1$, sign changes of the $\l_f(p)$ where $p$ goes over
primes have been considered by Ram Murty \cite{Mur}. Kohnen and Sengupta
\cite{KohSen} have shown that the first sign of $\l_f(n)$ happens for some
$n$ with
\bea\label{KohSen}
n\ll kN\exp\le(c\sqrt{\f{\log N}{\log\log(3N)}}\ri) (\log k)^{27}, \quad (n,N)=1,
\eea
where $c>2$ is a constant and the $\ll$-constant is absolute. Note that it is
natural to assume that $(n,N)=1$, since the eigenvalues $\l_f(p)$ with $p|N$
are explicitly known by the Atkin-Lehner theory.

Recently, Iwaniec, Kohnen, and Sengupta \cite{IwaKohSen} established the
following result.

\noindent
{\bf Proposition 2.7.} \text(Iwaniec-Kohnen-Sengupta \cite{IwaKohSen}).
{\it Let $f$ be a normalized Hecke eigenform of integral weight $k$ and
level $N$ that is a new form. Then there is some $n$ satisfying
\bea\label{IwaKohSen}
n\ll (k^2N)^{29/60}, \quad (n,N)=1,
\eea
such that $\l_f(n)<0$.}

The convexity bound
\bea\label{convL12f}
L(1/2+it,f)\ll (k^2N)^{1/4+\ve}
\eea
gives, instead of (\ref{IwaKohSen}), the weaker bound
\bea\label{consCONV}
n\ll (k^2N)^{1/2+\ve}, \quad (n,N)=1.
\eea
The uniform subconvexity bound
\bea\label{simSCB}
L(1/2+it,f)\ll (k^2N)^{29/120}
\eea
would prove Proposition 2.7, but no result of this quality is known.
The best known uniform subconvexity bound like (\ref{simSCB}) is
due to Michel and Venkatesh \cite{MicVen}, which states that
\bea\label{simSCB/MV/++}
L(1/2+it,f)\ll (k^2N)^{1/4-\d},
\eea
where $\d$ is some positive constant not specified.
Iwaniec, Kohnen, and Sengupta \cite{IwaKohSen}
manage to establish (\ref{IwaKohSen}) without appealing to (\ref{simSCB});
the key steps in \cite{IwaKohSen} are the following:
\begin{itemize}
  \item The identity (\ref{lf(p)^2=}), i.e.
$$
\l_f(p)^2=\l_f(p^2)+1, \quad (p,N)=1;
$$
  \item The Ramanujan conjecture (\ref{Ramun/Del}) proved by Deligne, i.e.
$$
|\l_f(p)|\leq 2, \quad (p,N)=1;
$$
  \item Sieve methods.
\end{itemize}
Of course, the proof of Proposition 2.7 is much more involved. In particular, to carry
out the sieves, one still needs the identity $\l_f(p)^2=\l_f(p^2)+1$ several times, and
needs the fact $|\l_f(p)|\leq 2$ in a more crucial way. We will not get into these
details, but just would like to point out that, the approach does not work for those
automorphic forms $f$, whose $\l_f(n)$ do not satisfy the above two properties.

\chapter{A Linnik-type problem for Maass forms}

\section{The spectral theory of Maass forms}
In this section, we
introduce the notion and basic facts from the theory of Maass forms of weight $k=0$ in the context of the Hecke
congruence subgroup $\G_0(N)$. Philosophically, there is no essential difference from the theory of classical
modular forms, except for the existence of a continuous spectrum in the space of Maass forms. A good monograph
on this topic is Iwaniec \cite{Iwa2}. But one has to admit that some mature methods, which are quite useful in
the case of holomorphic forms, do not work in the current situation. Linnik-type problem for Maass forms is such
an example, as will be explained in this chapter.

\subsection{The spectral decomposition: preliminary}
A function $f:\mathbb H\to \mathbb C$ is said to be
automorphic with respect to $\G_0(N)$ if
$$
f(\g z)=f(z), \quad \mbox{for all } \g\in \G_0(N).
$$
Therefore, $f$ lives on
$\G_0(N)\backslash\mathbb H$. We denote the space of such functions by
${\mathcal A}(\G_0(N)\backslash\mathbb H)$.
Our objective is to extend automorphic functions into automorphic forms subject to
suitable growth condition. The main results hold in the Hilbert space
$$
{\mathcal L}(\G_0(N)\backslash\mathbb H)=\{f\in {\mathcal A}(\G_0(N)\backslash\mathbb H):\|f\|<\infty\}
$$
with respect to the inner product
$$
\langle f,g \rangle =
\int_{\G_0(N)\backslash\mathbb H} f(z)\bar{g}(z)\f{{\rm d} x\,{\rm d} y}{y^2}.
$$
Recall that the standard Laplace operator on
the complex plane $\mathbb C $ is defined by
$$
\Delta^e=\frac{\partial^2}{\partial x^2}+\frac{\partial^2}{\partial y^2};
$$
but on the upper half-plane $\mathbb H$, we should use the non-Euclidean Laplace operator
$$
\Delta=-y^2 \bigg(\frac{\partial^2}{\partial x^2}+
\frac{\partial^2}{\partial y^2}\bigg).
$$
This non-Euclidean Laplace operator acts in the dense subspace
of smooth functions in $\mathcal L(\G_0(N)\backslash \mathbb H)$ such that
$f$ and $\Delta f$ are both bounded, i.e.
$$
{\mathcal D}(\G_0(N)\backslash\mathbb H)=\{f\in {\mathcal A}(\G_0(N)\backslash\mathbb H):
f, \D f \mbox{ smooth and bounded}\}.
$$
It is proved that ${\mathcal D}(\G_0(N)\backslash\mathbb H)$ is dense in ${\mathcal L}(\G_0(N)\backslash\mathbb H)$,
and $\D$ is positive semi-definite and symmetric on ${\mathcal D}(\G_0(N)\backslash\mathbb H)$. By Friedreich's
theorem in functional analysis, $\D$ has a unique self-adjoint extension to 
${\mathcal L}(\G_0(N)\backslash\mathbb H)$.

\noindent
{\bf Lemma 3.1.} {\it
\text{\rm (i)} Let $\L=s(1-s)$ be the eigenvalue of an
eigenfunction $f\in {\mathcal D}(\G_0(N)\backslash\mathbb H)$. Then $\L$ is real and non-negative, i.e.
either $s=1/2+it$ with $t\in \mathbb R$, or $0<s<1$.

\text{\rm (ii)} On ${\mathcal L}(\G_0(N)\backslash\mathbb H)$, the non-Euclidean Laplace operator $\D$ is positive
semi-definite and self-adjoint. }

With the above self-adjoint extension, one can show that
the non-Euclidean Laplace operator $\D$ has the spectral decomposition
$$
\mathcal L(\G_0(N)\backslash \mathbb H)=\mathbb C
\oplus \mathcal C(\G_0(N)\backslash \mathbb H)
\oplus \mathcal E(\G_0(N)\backslash \mathbb H).
$$
Here $\mathbb C$ is the space of constant functions,
$\mathcal C(\G_0(N)\backslash \mathbb H)$ the space of
cusp forms, and $\mathcal E(\G_0(N)\backslash \mathbb H)$
the space spanned by incomplete Eisenstein series.

\subsection{The discrete spectrum}
The structure of the space $\mathcal C(\G_0(N)\backslash \mathbb H)$,
the space of cusp forms, is characterized by the following result.

\noindent
{\bf Lemma 3.2.} {\it The automorphic Laplacian
$\D$ has a purely point spectrum on ${\mathcal C}(\G_0(N)\backslash\mathbb H)$, i.e.
the space ${\mathcal C}(\G_0(N)\backslash\mathbb H)$ is spanned by cusp forms. The
eigenvalues are
$$
0=\L_0<\L_1\leq \L_2\leq \ldots\to \infty,
$$
and the eigenspaces have finite dimension. For any complete orthonormal system
of cusp forms $\{u_j\}_{j=1}^\infty$, every $f\in {\mathcal C}(\G_0(N)\backslash\mathbb H)$
has the expansion
$$
f(z)=\sum_{j=1}^\infty \langle f,u_j\rangle u_j(z),
$$
converging in the norm topology. If
$f\in {\mathcal C}(\G_0(N)\backslash\mathbb H)\cap {\mathcal D}(\G_0(N)\backslash\mathbb H)$,
then the series converges absolutely and uniformly on compacta. }

Let
$$
\mathcal U=\{u_j\}_{j=1}^\infty
$$
be an orthonormal basis of $\mathcal C(\G_0(N)\backslash \mathbb H)$
which are eigenfunctions of $\D$, say
$$
\D u_j = \L_j u_j,
$$
with
\bea\label{lj=sj(1-sj)}
\L_j=s_j(1-s_j)=\f14+\nu_j^2, \quad s_j = \f12+i\nu_j.
\eea
Note that here the $\nu_j$ in (\ref{lj=sj(1-sj)}) is not necessarily real.
Any $u_j$ has the Fourier expansion
\bea\label{Four/uj}
u_j(z)=\sum_{n\not=0}\rho_j(n)W_{s_j}(nz),
\eea
where $W_s(z)$ is the Whittaker function given by
$$
W_s(z)=2|y|^{1/2}K_{s-1/2}(2\pi |y|)e(x),
$$
and $K_s(y)$ is the $K$-Bessel function. Note that
$$
W_s(z)\sim e(z), \quad y\to\infty.
$$
The automorphic forms $u_j(z)$ are called Maass
cusp forms. Sometimes, we write $f$ for Maass
cusp forms with Laplace eigenvalue
$$
\L=s(1-s)=\f14+\nu_f^2,
$$
and in this case, the Fourier expansion of $f$ takes the form
\bea\label{Four/uj/f}
f(z)=\sum_{n\not=0}\rho_f(n)W_{s}(nz).
\eea
Compare this with (\ref{Four/uj}).

\subsection{Antiholomorphic involution}
Let $\iota:\mathbb H\to \mathbb H$ be the antiholomorphic involution
$$
\iota(x+iy)=-x+iy.
$$
If $f$ is an eigenfunction of $\D$, and
\bea\label{Four/f}
f(z)=\sum_{n\not=0}\rho_f(n)W_{s}(nz),
\eea
then $f\circ\iota$ is an eigenfunction
with the same eigenvalue. Since $\iota^2=1$, its eigenvalues are $\pm 1$.
We may therefore diagonalize the Maass cusp forms with respect to
$\iota$. If $f\circ\iota=f$, we call $f$ even. In this case
$$
\rho_f(n)=\rho_f(-n).
$$
If $f\circ\iota=-f$, then we call $f$ odd, and we have
$$
\rho_f(n)=-\rho_f(-n).
$$

\subsection{The continuous spectrum}
On the other hand, in the space
${\mathcal E}(\G_0(N)\backslash\mathbb H)$, the spectrum
turns out to be continuous. The spectral resolution of $\D$ in ${\mathcal E}(\G_0(N)\backslash\mathbb H)$ follows from
the analytic continuation of the Eisenstein series.
The eigenpacket of the continuous spectrum consists
of the Eisenstein series $E_\mathfrak a(z,s)$ on the line $\s=1/2$ (analytically continued).
These are defined for every cusp $\mathfrak a$ by
$$
E_\mathfrak a(z,s)=\sum_{\g\in \G_\mathfrak a\backslash \G_0(N)}(\Im \s_{\mathfrak a}^{-1}\g z)^s
$$
if $\s>1$, and by analytic continuation for all $s\in \mathbb C$. Here $\G_\mathfrak a$ is the
stability group of $\mathfrak a$ and $\mathfrak a\in SL_2(\mathbb R)$ is such that
$$
\s_\mathfrak a \infty =\mathfrak a, \quad \s_{\mathfrak a}^{-1}\G\s_{\mathfrak a}
=\G_\infty.
$$
The scaling matrix $\s_{\mathfrak a}$ of $\mathfrak a$ is only determined
up to a translation from the right; however the Eisenstein series
does not depend on the choice of $\s_{\mathfrak a}$, not even
on the choice of a cusp in the equivalent class. The Fourier expansion of
$E_\mathfrak a(z,s)$ is similar to that of a cusp form; precisely,
$$
E_\mathfrak a(z,s)=\vp_\mathfrak a y^s + \vp_\mathfrak a (s)y^{1-s}
+\sum_{n\not=0}\vp_\mathfrak a (n,s)W_s(nz),
$$
where $\vp_\mathfrak a=1$ if $\mathfrak a \sim\infty$, and
$\vp_\mathfrak a=0$ otherwise.

\noindent
{\bf Lemma 3.3.} {\it The space ${\mathcal E}(\G_0(N)\backslash\mathbb H)$
of incomplete Eisenstein series splits orthogonally into
$\D$-invariant subspaces
$$
{\mathcal E}(\G_0(N)\backslash\mathbb H)=
{\mathcal R}(\G_0(N)\backslash\mathbb H)\oplus_{\mathfrak a}
{\mathcal E}_{\mathfrak a}(\G_0(N)\backslash\mathbb H).
$$
The spectrum of $\Delta$ in ${\mathcal R}(\G_0(N)\backslash\mathbb H)$ is discrete; it consists
of a finite number of points $\L_j$ with
$$
\L_j\in [0,1/4).
$$
The spectrum of
$\D$ on ${\mathcal E}_{\mathfrak a}(\G_0(N)\backslash\mathbb H)$ is absolutely continunous;
it covers the segment
$$
[1/4,+\infty)
$$
uniformly with multiplicity $1$. Every $f\in {\mathcal E}(\G_0(N)\backslash\mathbb H)$ has the expansion
\bea
f(z)=\sum_{j}\langle f,u_j\rangle u_j(z)+\sum_{\mathfrak a}\f{1}{4\pi}\int_{-\infty}^{\infty}
\langle f,E_{\mathfrak a}(\cdot, 1/2+it)\rangle E_{\mathfrak a}(z,1/2+it){\rm d}t,
\eea
which converges in the norm topology. If
$f\in {\mathcal E}(\G_0(N)\backslash\mathbb H)\cap {\mathcal D}(\G_0(N)\backslash\mathbb H)$,
then the series converges pointwise absolutely and uniformly on compacta.
}

\subsection{The spectral decomposition: conclusion}
Combining Lemmas 3.2-3.3, one gets the spectral decomposition of the
whole space ${\mathcal L}(\G_0(N)\backslash\mathbb H)$,
\bna
f(z)
&=&\sum_{j=0}^\infty \langle f,u_j\rangle u_j(z)
+\sum_{j} \langle f,u_j\rangle u_j(z) \\
&& +\sum_{\mathfrak a}\f{1}{4\pi}\int_{-\infty}^{\infty}
\langle f,E_{\mathfrak a}(\cdot, 1/2+it)\rangle E_{\mathfrak a}(z,1/2+it)dt.
\ena
This structure is one of the basics for later arguments.

\section{Hecke theory for Maass forms}

\noindent
For $n\geq 1$, define
\bea\label{Tnfz=MAASS}
(T_n f)(z)=\f{1}{\sqrt{n}}\sum_{ad=n}
\sum_{b\,({\rm mod}\,d)} f\le(\f{az+b}{d}\ri);
\eea
nevertheless, only those $T_n$ with $(n,N)=1$ are interesting.
We first examine the action of $T_n$ on a Maass cusp form $u_j$.
For $u_j$ as in (\ref{lj=sj(1-sj)}) and (\ref{Four/uj}), write
\bea\label{fourierMAASS}
u_j(z)=\sum_{m\neq 0}\rho_j(m)W_{s_j}(mz).
\eea
Then one computes that
$$
(T_n u_j)(z)=\sum_{m\neq 0}t_n(m)W_{s_j}(mz),
$$
with
$$
t_n(m)=\sum_{d|(m,n)}\rho_j\le(\f{mn}{d^2}\ri).
$$
It follows that
\bna
T_m T_n=\sum_{d|(m,n)}T_{mn/d^2},
\ena
so that in particular $T_m$ and $T_n$ commute. Moreover, the Hecke operators
commute with the non-Euclidean Laplace operator $\D$. For every $n$,
$T_n$ is also self-adjoint on ${\mathcal L}(\G_0(N)\backslash \mathbb H)$, i.e.
$$
\langle T_n f, g\rangle=\langle f, T_n g\rangle, \quad (n,N)=1.
$$
Therefore, in the space ${\mathcal C}(\G_0(N)\backslash \mathbb H)$ of cusp forms,
an orthonormal basis $\{u_j\}_{j=1}^\infty$ can be chosen which consists of simultaneous
eigenfunctions for all $T_n$, i.e.
\bea\label{Tnuj=ljnuj}
T_n u_j(z)=\l_j(n)u_j(z), \quad j\geq 1, \ (n,N)=1,
\eea
where $\l_j(n)$ is the eigenvalue of $T_n$ for $u_j(z)$.
Up to a constant, $\l_j(n)$ and the Fourier coefficient $\rho_j(n)$ are equal. More precisely,
\bea\label{ljnrhoj1}
\l_j(n)\rho_j(1)= \rho_j(n), \quad \mbox{for all } (n,N)=1, \ j\geq 1.
\eea
Note that if $\rho_j(1)=0$, then all $\rho_j(n)=0$ for $(n,N)=1$.

\noindent
{\bf Lemma 3.4.} {\it Let
$
{\mathcal U}=\{u_j\}_{j=1}^\infty
$
be an orthonormal basis consisting of simultaneous
eigenfunctions for all $T_n$. Fix a $j$, and let
$\{\l_j(n)\}_{n=1}^\infty$ be the sequence of eigenvalues for all $T_n$ as in
(\ref{Tnuj=ljnuj}).

\vskip -2mm

\mbox{\rm (i)} The Hecke eigenvalues $\{\l_j(n)\}_{n=1}^\infty$ are real;

\vskip -2mm

\mbox{\rm (ii)}  The Hecke eigenvalues are multiplicative
in the following sense:
$$
\l_j(m)\l_j(n)=\sum_{d|(m,n)}\l_j\le(\f{mn}{d^2}\ri),
\quad (n,N)=1,
$$
and
\bea\label{lf(mp)/MAASS}
\l_j(m)\l_j(p)=\l_j(mp), \quad p|N.
\eea
It follows that
\bea\label{lf(p)^2=MAASS}
\l_j(p)^2=\l_j(p^2)+1, \quad (p,N)=1.
\eea}

As in the case of classical modular forms, we cannot deduce from (\ref{ljnrhoj1}) that
$\rho_j(n)\not=0$. Thus, we need to work with the new forms for the same reason.

\noindent
{\bf Lemma 3.5.} {\it If $u_j$ is a new form, then
(\ref{Tnuj=ljnuj}) holds for all $n$. The first coefficient
in the Fourier expansion (\ref{fourierMAASS}) does not vanish,
so one can normalize $u_j$ by setting $\rho_j(1)=1$. In this case,
$\rho_j(n)=\l_j(n)$ for all $n$, and hence
\bea\label{fourierMAASS/NEW}
u_j(z)=\sum_{n\not=0}^\infty \l_j(n)W_{s_j}(nz).
\eea
 }

The Eisenstein series $E_\infty (z,1/2+it)$ is an eigenfunction of all
the Hecke operators $T_n, (n,N)=1,$ with eigenvalues
$\eta_t(n),$
i.e.
\bea
T_n E_\infty(z,1/2+it)=\eta_t(n)E_\infty(z,1/2+it), \quad (n,N)=1, \ t \in \mathbb R,
\eea
where
\bea
\eta_t(n)=\sum_{ad=n}\le(\f{a}{d}\ri)^{it}.
\eea
We recall that $\eta_0(n)$ reduces to the classical divisor function $\tau(n)$.

\section{Automorphic $L$-functions for Maass forms}

\noindent
To a Maass
new form $f$ as in (\ref{Four/uj/f}) with Laplace eigenvalue $1/4+\nu^2$, we may attach an automorphic
$L$-function \bea\label{def/Lsf/MAASS} L(s,f)=\sum_{n=1}^\infty \f{\l_f(n)}{n^s}, \qquad \s>1, \eea as in \S
2.2. For the Maass case, the complete $L$-function is defined as \bea\label{def/compLsf/Maass}
\Phi(s,f)=\pi^{-s}\G\le(\f{s+\epsilon-1/2+\nu}{2}\ri) \G\le(\f{s+\epsilon+1/2-\nu}{2}\ri) L(s,f), \eea where
$\epsilon$ is the eigenvalue of $\iota$ introduced in \S3.1.3. The complete $L$-function satisfies the
functional equation 
$$
\Phi(s,f)=\ve_f N^{1/2-s}\Phi(1-s,\bar{f}), 
$$
where $\ve_f$ is a
complex number of modulus $1$. For any new form $f$, we have the following Euler product for $\Phi(s,f)$.

\noindent
{\bf Lemma 3.6.} {\it If $f$ is a new form, then the functional
equation takes the form
\bea\label{funeqn/NEW/Maass}
\Phi(s,f)=\ve_f N^{1/2-s}\Phi(1-s,f),
\eea
and, for $\s>1,$ the function $L(s,f)$ admits the Euler product
\bea\label{EuProd/MAASS}
L(s,f)=\prod_{p}\le(1-\f{\l_f(p)}{p^s}+\f{\chi_N^0(p)}{p^{2s}}\ri)^{-1},
\eea
where $\chi_N^0$ is the principal character modulo $N$.}

By Lemma 3.5, all eigenvalues $\l_f(n)$ of a new form $f$ are real.
This explains why on the right-hand side of (\ref{funeqn/NEW/Maass}) we have
$\Phi(1-s,f)$ instead of $\Phi(1-s,\bar f)$.
We may further factor the Hecke polynomial in (\ref{EuProd/MAASS}) into
$$
1-\f{\l_f(p)}{p^s}+\f{\chi_N^0(p)}{p^{2s}}=
\le(1-\f{\a_f(p)}{p^s}\ri)\le(1-\f{\b_f(p)}{p^s}\ri),
$$
where
$$
\a_f(p)+\b_f(p)=\l_f(p), \quad \a_f(p)\b_f(p)=\chi^0_N(p).
$$
The Generalized Ramanujan's Conjecture (GRC in brief) in this case asserts that
\bea\label{Ramun/Del/MAASS}
|\a_f(p)|=|\b_f(p)|=1, \quad (p,N)=1;
\eea
this is still open, and the strongest bound towards the above conjecture
is that of Kim and Sarnak \cite{KimSar}:
\bea
|\a_f(p)|\leq p^{7/64},
\qquad
 |\b_f(p)|\leq p^{7/64},
 \qquad (p,N)=1.
\eea
The symmetric square
$L$-function is defined by
\bea\label{defsym2/MAASS}
L(s,\text{\rm sym}^2f)=L(2s,\chi_N^0)\sum_{n=1}^\infty \f{\l_f(n^2)}{n^s}.
\eea
Now suppose $g$ is a Maass new form of level $N'$ and weight $k'$.
The Rankin-Selberg $L$-function of $f$ and $g$ is defined as
\bea\label{defRSL/MAASS}
L(s,f\otimes g)=L(2s,\chi_N^0\chi_{N'}^0)\sum_{n=1}^\infty \f{\l_f(n)\l_g(n)}{n^s},
\eea
if $[N,N']$ is square-free. For $\s>1$, the Euler product
of $L(s,f\otimes g)$ takes the form
\bea\label{RSLeuler/MAASS}
L(s,f\otimes g)
&=&\prod_{p}\le(1-\f{\a_f(p)\a_g(p)}{p^{s}}\ri)^{-1}
\le(1-\f{\a_f(p)\b_g(p)}{p^{s}}\ri)^{-1}\nonumber\\
&&\times \le(1-\f{\b_f(p)\a_g(p)}{p^{s}}\ri)^{-1}
\le(1-\f{\b_f(p)\b_g(p)}{p^{s}}\ri)^{-1}.
\eea

The following subconvexity bound is a new result of
Michel and Venkatesh \cite{MicVen}.

\noindent
{\bf Lemma 3.7.} \text{(Michel-Venkatesh \cite{MicVen}).}
{\it Let $f$ be a non-zero
Maass new form on $\G_0(N)$ with Laplace eigenvalue $1/4+\nu^2$.
Then
\bea\label{conv/MAASS}
L(1/2+it, f)\ll \{(1+|t+\nu|)^2N\}^{1/4-\d/2},
\eea
where $\d$ is a positive absolute constant. }

Subconvexity bounds for any one of the three aspects $\nu, N,$ or $t$
have been studied extensively in the literature, but uniform subconvexity bound is only
known of the shape (\ref{conv/MAASS}), where $\d>0$ is not specified.
See Michel \cite{Mic} for a survey in this direction, and Michel and Venkatesh \cite{MicVen}
for recent developments.

\section{Infinite sign changes of Fourier coefficients of Maass forms}

\noindent
It is pointed out at the end of \cite{KnoKohPri} that
similar results hold for Maass forms. Therefore, we have the following
general result.

\noindent
{\bf Proposition 3.8.} \text{(Knopp-Kohnen-Pribitkin \cite{KnoKohPri}).}
{\it Let $f$ be a non-zero
Maass cusp form on $\G_0(N)$ with Fourier expansion (\ref{Four/f}),
and suppose that its Fourier coefficients $\rho_f(n)$ are real for all
$n\geq 1$. Then  the sequence
$$
\{\rho_f(n)\}_{n=1}^\infty
$$
has infinitely many sign changes, i.e. there are
infinitely many $n$ such that $\rho_f(n)>0$, and there are infinitely
many $n$ such that $\rho_f(n)<0$.
}

Since all eigenvalues $\l_f(n)$ of a new form $f$
are real, we arrive at the following corollary.

\noindent
{\bf Corollary 3.9.} {\it Let $f$ be a normalized Maass
eigenform that is a new form of level $N$ on $\G_0(N)$.
Then  the sequence
$
\{\l_f(n)\}_{n=1}^\infty
$
has infinitely many sign changes, i.e. there are
infinitely many $n$ such that $\l_f(n)>0$, and there are infinitely
many $n$ such that $\l_f(n)<0$.
}

For Maass forms, one may also ask the more precise question
how long is the sequence of $\{\l_f(n)\}_{n=1}^\infty$
that keep the same sign. Like in (\ref{def/C+-}), one may also
introduce
\bea\label{def/C+-/MAASS}
{\mathscr N}_f^+(x)=\sum_{\substack{n\leq x, \, (n,N)=1\\ \l_f(n)>0}} 1,
\eea
and define ${\mathscr N}_f^-(x)$ similarly by replacing the condition $\l_f(n)>0$
under the summation by $\l_f(n)<0$. It is possible to establish results
similar to those in Kohnen-Lau-Shparlinski \cite{KohLauShp},
in Wu \cite{Wu}, or even in Lau-Wu \cite{LauWu}. But this will
carry us too far, and we prefer to do it elsewhere at a later stage.

\section{A Linnik-type problem for Maass forms}

\noindent
As in the case of holomorphic eigenforms, one may also formulate
Linnik's problem for Maass eigenforms. That is:

{\it For a Maass eigenform $f$, is it possible to obtain a bound on the
first sign change of $\l_f(n)$, say,
in terms of $N$ and the Laplace eigenvalue of $f$? }

There seems no result in this direction. Our result in the following
produces one.

\noindent
{\bf Theorem 3.10.}
{\it Let $f$ be a normalized Maass new form of
level $N$ and Laplace eigenvalue $1/4+\nu^2$.
Then there is some $n$ satisfying
\bea\label{n<</MAASS}
n\ll ((3+|\nu|)^2N)^{1/2-\delta}, \quad (n,N)=1,
\eea
such that $\l_f(n)<0$, where $\delta$ is a positive absolute constant. }

We need Perron's formula in the following form, the proof of which can be
found in standard text books on analytic number theory.

Let $\s_a$ be the abscissa
of absolute convergence for
the Dirichlet series
\bea\label{1/F(s)}
F(s)=\sum_{n=1}^\infty \f{a_n}{n^s},
\eea
where $\{a_n\}_{n=1}^{\infty}$ is a sequence of
complex numbers, and $s=\s+it\in \mathbb C$ a complex variable.
Perron's formula expresses a partial sum of the coefficients
$a_n$ in terms of $F(s)$.

\noindent
{\bf Lemma 3.11.} \text{\rm (Perron's formula).}
{\it Define
\bea\label{1/AB}
A(x)=\max_{x/2<n\leq 3x/2} |a_n|, \quad
B(\s)=\sum_{n=1}^{\infty} \f{|a_{n}|}{n^{\s}}
\eea
for $\s>\s_a$. Let $\ell$ be a non-negative integer,
$x\geq 2$, and $\|x\|$ denote the distance between $x$ and the nearest integer.
Then, for $b>\s_a$ and $T\geq 2$,
\bea\label{1/per}
\sum_{n\leq x}a_n\bigg(\log \f{x}{n}\bigg)^{\ell}
&=&\f{1}{2\pi i}\int_{b-iT}^{b+iT}F(s)\f{x^s}{s^{\ell+1}}{\rm d} s
+O\le(\f{xA(x)\log^{\ell+1} x}{T}\ri)\nonumber\\
&&+O\le(\f{x^{b}B(b)\log^\ell x}{T}\ri) \nonumber\\
&&+O\le\{A(x)\min\le(1,\f{x}{T\|x\|}\ri)\log^\ell x\ri\}.
\eea
In particular, for $b>\s_a$,
\bea\label{1/perSIM}
\sum_{n\leq x}a_n\bigg(\log \f{x}{n}\bigg)^{\ell}
=\f{1}{2\pi i}\int_{b-i\infty}^{b+i\infty}F(s)\f{x^s}{s^{\ell+1}}{\rm d} s.
\eea
}

Now we are in a position to establish Theorem 3.10.

\medskip
\noindent
{\it Proof of Theorem 3.10.} The idea is to consider the sum
$$
S(x):=\sum_{\substack{n\leq x\\ (n,N)=1}} \l_f(n)\log \f{x}{n},
$$
assuming that
\bea\label{assum/MAASS}
\l_f(n)\geq 0 \quad \text{\rm for} \ n\leq x \ \text{\rm and }  (n,N)=1.
\eea
The desired result will follow from upper and lower bound estimates
for $S(x)$.

To get an upper bound for $S(x)$, we apply Perron's formula (\ref{1/perSIM}) with $\ell=1$
to the Dirichlet series (\ref{def/Lsf/MAASS}), getting
$$
\sum_{n\leq x} \l_f(n)\log \f{x}{n}
=\f{1}{2\pi i}\int_{2-i\infty}^{2+i\infty}L(s,f)\f{x^s}{s^2}{\rm d} s.
$$
Moving the contour to the vertical line $\s=1/2$, where we apply the
Michel-Venkatesh bound (\ref{conv/MAASS}) for $L(s,f)$, we obtain
\bna
\sum_{n\leq x} \l_f(n)\log \f{x}{n}
&=&\f{1}{2\pi i}\int_{1/2-i\infty}^{1/2+i\infty}L(s,f)\f{x^s}{s^2}{\rm d} s
\\
&\ll& \int_{-\infty}^{\infty}\{(1+|t+\nu|)^2N\}^{1/4-\d/2}\f{x^{1/2}}{|t|^2+1}{\rm d} t
\\\noalign{\vskip 2mm}
&\ll& ((3+|\nu|)^2N)^{1/4-\d/2}x^{1/2}.
\ena
To recover an estimate for $S(x)$ from the above result, we introduce the condition
$(n,N)=1$ by means of the M\"obius inversion formula, which gives
$$
S(x)=\sum_{d|N}\mu(d)\sum_{dm\leq x}\l_f(dm)\log\f{x}{dm}.
$$
Since $d|N$, we may apply the mutiplicativity property (\ref{lf(mp)/MAASS}), which states
that in the current situation
$$
\l_f(dm)=\l_f(d)\l_f(m).
$$
It follows that
\bna
S(x)
&=&\sum_{d|N}\mu(d)\l_f(d)\sum_{dm\leq x}\l_f(m)\log\f{x}{dm} \\
&\ll &\sum_{d|N}|\mu(d)\l_f(d)| \bigg|\sum_{m\leq x/d}\l_f(m)\log\f{x/d}{m}\bigg| \\
&\ll & ((3+|\nu|)^2N)^{1/4-\d}x^{1/2} \sum_{d|N}\f{|\l_f(d)|}{d^{1/2}}.
\ena
The Kim-Sarnak bound states that $\l_f(d)\ll d^{7/64+\ve}$, and therefore,
$$
\sum_{d|N}\f{|\l_f(d)|}{d^{1/2}}\ll \tau(N)\ll N^\ve;
$$
here we note that the trivial bound $\l_f(d)\ll d^{1/2+\ve}$ works equally well. Consequently,
we conclude that
\bea\label{uppMAASS}
S(x)\ll ((3+|\nu|)^2N)^{1/4-\d/2}x^{1/2}.
\eea

To get a lower bound
for $S(x)$ under the assumption (\ref{assum/MAASS}),
we first get rid of the weight $\log (x/n)$ in a simple way:
$$
S(x)\gg \sum_{\substack{n\leq x/2\\ (n,N)=1}} \l_f(n).
$$
We now restrict the summation to $n=pq$, where $p$ and $q$ are primes
satisfying
$$
p\leq \sqrt{x/2}, \quad
q\leq \sqrt{x/2}, \quad
(p,N)=1, \quad
(q,N)=1,
$$
and use the formulae
\bna
\left\{
\begin{array}{lll}
\l_f(pq)=\l_f(p)\l_f(q) & \text{if } p\not=q, \ (p,N)=1, \ (q,N)=1,
\\\noalign{\vskip 3mm}
\l_f(p^2)=\l_f(p)^2-1  & \text{if } p=q, \ (p,N)=1.
\end{array}
\right.
\ena
We get
\bna
S(x)
&\gg& \sum_{\substack{p\leq \sqrt{x/2}\\ (p,N)=1}}
\sum_{\substack{q\leq \sqrt{x/2}\\ (q,N)=1}}
\l_f(pq) \\
&=&\bigg\{\sum_{\substack{p\leq \sqrt{x/2}\\ (p,N)=1}}\l_f(p)\bigg\}^2
-\sum_{\substack{p\leq \sqrt{x/2}\\ (p,N)=1}}1.
\ena
Recalling the assumption (\ref{assum/MAASS}), we have
$\l_f(p^2)\geq 0$ for $p\leq \sqrt{x/2}$ and $(p,N)=1$, and therefore,
$$
\l_f(p)^2=\l_f(p^2)+1\geq 1,
$$
that is $\l_f(p)\geq 1.$ It follows from this and the prime number theorem that
\bea\label{lowMAASS}
S(x)&\geq & \bigg\{\sum_{p\leq \sqrt{x/2}}1\bigg\}^2-\bigg\{\sum_{p\leq \sqrt{x/2}}1 \bigg\}\nonumber\\
&\gg & \f{x}{\log^2x}.
\eea

Comparing (\ref{lowMAASS}) with (\ref{uppMAASS}), we get
$$
\f{x}{\log^2x}\ll S(x)\ll ((3+|\nu|)^2N)^{1/4-\d/2}x^{1/2},
$$
that is
$$
x\ll ((3+|\nu|)^2N)^{1/2-\d+\ve}.
$$
This proves the theorem.
\hfill $\square$

\chapter{A Linnik-type problem for automorphic $L$-functions}

\section{Automorphic $L$-functions: concepts and properties}

\noindent
To each irreducible unitary cuspidal representation
$\pi=\otimes \pi_p$ of $GL_m({\mathbb A}_{\mathbb Q})$, one can attach a global $L$-function
$L(s,\pi)$, as in Godement and Jacquet \cite{GodJac}, and Jacquet and
Shalika \cite{JacSha1}. For $\s=\Re s>1$, $L(s,\pi)$ is defined by
products of local
factors
\bea\label{Lspi}
L(s,\pi)
=\prod_{p<\infty} L_p (s, \pi_p),
\eea
where
\bea\label{Lpspip}
L_p(s, \pi_p )=\prod_{j=1}^m
\le(1-\f{\alpha_\pi(p,j)}{p^{s}}\ri)^{-1};
\eea
the complete $L$-function $\Phi(s,\pi)$ is defined by
\bea\label{PHIspi}
\Phi(s,\pi)=L_\infty(s,\pi_\infty)L(s,\pi),
\eea
where
\bea\label{L8spi8}
L_\infty(s, \pi_\infty)=
\prod_{j=1}^m
\Gamma_{\mathbb R}(s+\mu_{\pi}(j))
\eea
is the Archimedean local factor. Here
\bea\label{GammaR2}
\Gamma_{\mathbb R}(s) =\pi^{-s/2}\Gamma \le(\f{s}{2}\ri),
\eea
and $\{\alpha_\pi(p,j)\}_{j=1}^m$ and $\{\mu_{\pi} (j)\}_{j=1}^m $
are complex numbers associated with $\pi_p$ and $\pi_\infty$,
respectively, according to the Langlands correspondence.
The case $m=1$ is classical; for $m\geq 2$, $\Phi(s,\pi)$ is
entire and satisfies a functional equation.

We review briefly some properties of the automorphic $L$-functions
$L(s,\pi)$ and $\Phi(s,\pi)$, which we will use for our proofs.

\noindent
{\bf Lemma 4.1.} (Jacquet-Shalika \cite{JacSha1}).
{\it The Euler product for $L(s,\pi)$ in
(\ref{Lspi}) converges absolutely for $\s>1$.
}

Thus, in the half-plane $\s>1$, we may write
\bea\label{Lspi=Dir}
L(s,\pi)=\sum_{n=1}^\infty \f{\l_\pi(n)}{n^s},
\eea
where
\bea\label{lpin=}
\l_\pi(n)=\prod_{p^\nu\|n}
\bigg\{\sum_{\nu_1+\cdots +\nu_m=\nu} \a_\pi(p,1)^{\nu_1}\cdots \a_\pi(p,m)^{\nu_m} \bigg\}.
\eea
In particular,
\bea\label{lpin=PAR}
\l_\pi(1)=1, \quad \l_\pi(p)=\a_\pi(p,1)+\cdots +\a_\pi(p,m).
\eea

\noindent
{\bf Lemma 4.2.} (Shahidi \cite{Sha1}, \cite{Sha2}, \cite{Sha3}, and \cite{Sha4}).
{\it The complete $L$-function
$\Phi(s,\pi)$ has an analytic continuation to
the whole complex plane and satisfies the functional equation
$$
\Phi(s,\pi) =\ve(s,\pi)
\Phi(1-s,\tilde{\pi})
$$
with
$$
\ve(s,\pi)=
\ve_\pi N_{\pi}^{1/2-s},
$$
where $N_{\pi}\geq 1$ is an integer called the arithmetic conductor of $\pi$, $\ve_\pi$
is the root number satisfying $|\ve_\pi|=1$,
and $\tilde{\pi}$ is the representation contragredient to $\pi$.}

If $\tilde{\pi}$ is the representation contragredient to $\pi$, then we have
\bea\label{4/for9}
\{\a_{\tilde\pi}(p,j)\}_{j=1}^m
=\{\overline{\a_\pi(p,j)}\}_{j=1}^m
\eea
and
\bea
\{\mu_{\tilde\pi}(j)\}_{j=1}^m
=\{\overline{\mu_\pi(j)}\}_{j=1}^m.
\eea
It follows from these and (\ref{lpin=}) that
\bea\label{l=lover}
\l_{\tilde\pi}(n)=\overline{\l_{\pi}(n)}.
\eea
Therefore, if $\pi$ is self-contragredient, i.e.
$\pi=\tilde{\pi}$, then (\ref{l=lover}) states that
\bea\label{reall}
\l_{\pi}(n)=\overline{\l_{\pi}(n)},
\eea
which means that $\l_{\pi}(n)$ is real.

\noindent
{\bf Lemma 4.3.} (Godement-Jacquet \cite{GodJac},
and Jacquet-Shalika \cite{JacSha1}).
{\it The function
$\Phi(s,\pi)$ is entire, and bounded in vertical strips with finite
width.}

\noindent
{\bf Lemma 4.4.} (Gelbart-Shahidi \cite{GelSha}).
{\it The function
$\Phi(s,\pi)$ is of order one.}

\noindent
{\bf Lemma 4.5.}  (Jacquet-Shalika \cite{JacSha1},
and Shahidi \cite{Sha1}).
{\it The function $\Phi(s,\pi)$ and
$L(s,\pi)$ are non-zero in the half-plane $\s\ge 1$. }

Iwaniec and Sarnak \cite{IwaSar} introduced the analytic
conductor of $\pi$. It is a function over the reals given by
\bea\label{AnalCod}
Q_\pi(t)=N_\pi \prod_{j=1}^m (3+|t+\mu_\pi(j)|),
\eea
which puts together all the important parameters for $\pi$. The
quantity
\bea\label{Cod}
Q_\pi:=Q_\pi(0)=N_\pi \prod_{j=1}^m (3+|\mu_\pi(j)|)
\eea
is also important, and it is named as the conductor of $\pi$.

\goodbreak

The next lemma is about the distribution of zeros of the
function $L(s,\pi)$.

\noindent
{\bf Lemma 4.6.} {\it All the non-trivial zeros of $\Phi(s,\pi)$ are in the critical
strip $0\leq \s \leq 1$. Let $N(T,\pi)$ be the number of nontrivial zeros within the rectangular
$$
0\leq \s \leq 1, \ |t|\leq T.
$$
Then
$$
N(T,\pi)\ll T\log (Q_{\pi}T),
$$
and
$$
N(T+1,\pi)-N(T,\pi)\ll \log (Q_{\pi}T).
$$
}

For proof of this, one is referred to Liu and Ye \cite{LiuYe1}, Lemma 4.3,
or Iwaniec and Kowalski \cite{IwaKow}, Theorem 5.8.

\section{Three conjectures in the theory of automorphic $L$-functions}

It is said in the previous section that
all the non-trivial zeros of $\Phi(s,\pi)$ are in the critical
strip $0\leq \s \leq 1$, while GRH for $L(s,\pi)$
predicts that they should actually lie on the
vertical line $\s=1/2$.

\noindent
{\bf Generalized Riemann Hypothesis.} {\it All the zeros
of $\Phi(s,\pi)$ have their real parts equal to $1/2$. }

Upper bounds for $L(s,\pi)$ on the critical line $\s=1/2$ is of great importance,
and the most optimistic conjecture in this direction can be stated as follows
in terms of the
analytic conductor defined in (\ref{AnalCod}).

\noindent
{\bf Generalized Lindel\"of Hypothesis.}
{\it The estimate
$$
L\le(\f12+it,\pi\ri)\ll Q_\pi(t)^\ve
$$
is true for arbitrary $\ve>0$.}

The following result is unconditional. Its proof is based on the
fact that the Rankin-Selberg $L$-function $L(s,\pi\otimes\pi')$
exists, where $\pi'$ is another irreducible unitary cuspidal representation.
It may happen that $\pi=\pi'$.

\noindent
{\bf Lemma 4.7.} (Harcos \cite{Har}).
{\it Let $\ve>0$ be arbitrary, and $0<\sigma<1$. Then we have the upper bound estimate
\bea\label{conv/AUTO}
L(\s+it,\pi)\ll_\ve Q_\pi(t)^{\f{1-\s}{2}+\ve}.
\eea
}

Taking $\s=1/2$, Lemma 4.7 gives
$$
L\le(\f12+it,\pi\ri)\ll_\ve Q_\pi(t)^{\f{1}{4}+\ve}.
$$
This is called the convexity bound of $L(s,\pi)$, and it should be
emphasized that it is uniform in all parameters. Subconvexity bounds 
have been established for some $L(s,\pi)$ for some aspects, only when
$m=1,2,3,4,8$; moreover, when $m\geq 2$, the existing subconvexity
bounds are not uniform in all parameters, except the recent uniform result
of Michel and Venkatesh \cite{MicVen} for $GL_2$. See Michel \cite{Mic} for a
comprehensive survey in this direction.

\medskip

Good bounds for the local parameters
$$
\{\a_\pi(p,j)\}_{j=1}^m, \qquad \{\mu_\pi(p,j)\}_{j=1}^m
$$
are of fundamental importance for the study of automorphic $L$-functions. By the Rankin-Selberg
method, one shows that, for all $p$,
\bea\label{GRCtri}
|\a_\pi(p,j)|\leq
p^{1/2}, \quad \Re \mu_\pi(j)\leq \f12;
\eea
moreover, for any unramified place,
\bea\label{GRCtri+}
p^{-1/2}\leq |\a_\pi(p,j)|\leq
p^{1/2}, \quad |\Re \mu_\pi(j)|\leq \f12.
\eea
The bounds (\ref{GRCtri}) and (\ref{GRCtri+}) are called trivial
bounds and are hence of little use.
The Generalized Ramanujan Conjecture (GRC in brief) asserts that
the $1/2$ in (\ref{GRCtri+}) can be reduced to $0$.

\noindent
{\bf Generalized Ramanujan Conjecture.} {\it With $\a_\pi(p,j)$
and $\mu_\pi(j)$ defined as above,
$$
\left\{\begin{array}{rl}
|\a_\pi(p,j)|=1 & \mbox{if $\pi$ is unramified at $p$},
\\\noalign{\vskip 3mm}
|\Re \mu_\pi(j)|=0 & \mbox{if $\pi$ is unramified
at $\infty$}.
\end{array}
\right.
$$}

The following lemma gives bounds toward the GRC.

\noindent
{\bf Lemma 4.8.} (Luo-Rudnick-Sarnak \cite{LuoRudSar}).
{\it There is a constant $0\leq \th<1/2$, such that
$$
\left\{
\begin{array}{rl}
|\a_\pi(p,j)|\leq p^\th & \mbox{if $\pi$ is unramified at $p$},
\\\noalign{\vskip 3mm}
|\Re \mu_\pi(j)|\leq \th & \mbox{if $\pi$ is unramified
at $\infty$}.
\end{array}
\right.
$$
Actually,
\bea\label{bdGRC}
\th=\f12-\f{1}{m^2+1}
\eea
is acceptable. }

This $\th$ will be used in the next two chapters.

\section{Hecke $L$-functions as automorphic $L$-functions}

\noindent
It should be pointed out that the Hecke $L$-functions defined in
(\ref{def/Lsf}) and (\ref{def/Lsf/MAASS}) are special examples of automorphic $L$-functions.
If $\pi$ corresponds to holomorphic new form $f$ with weight $k$ and level $N$,
then the conductor is
$$
Q_\pi\asymp k^2N.
$$
If $\pi$ corresponds to Maass new form $f$ with Laplace eigenvalue $1/4+\nu^2$ and level $N$, then
$$
Q_\pi\asymp (3+|\nu|)^2N.
$$
In view of these, it is easy to re-state Theorem 3.10
in terms of their conductors $Q_\pi$.

\section{Rankin-Selberg $L$-functions}
Let $\pi$ and $\pi'$ be two
irreducible unitary cuspidal representations
for $GL_m(\mathbb A_\mathbb Q)$ and
$GL_{m'}(\mathbb A_\mathbb Q)$, respectively. The theory for the
Rankin-Selberg type $L$-functions $L(s,\pi\otimes \pi')$
was initiated by Rankin \cite{Ran} and Selberg
\cite{Sel2} in the
case of classical modular forms. For general
automorphic representations, the corresponding theory
was initiated and developed in several papers
by Jacquet, Pisteski-Shapiro, and Shalika
\cite{Jac} \cite{JacSha1} \cite{JacSha2},
and completed in works of Shahidi \cite{Sha1}
\cite{Sha2} \cite{Sha3} \cite{Sha4},
Moeglin and Waldspurger \cite{MoeWal},
and Gelbart and Shahidi \cite{GelSha}.
Let $\pi$ and $\pi'$ be as above. When $\s>1$,
\bea
L(s,\pi\otimes\pi')=
\prod_{p<\infty}L_p(s,\pi_p\otimes\pi'_p)
\eea
with
$$
L_p(s,\pi_p\otimes\pi'_p)=
\prod_{j=1}^{mm'}\le(1-\f{\a_{\pi\otimes\pi'}(p,j)}{p^s}\ri)^{-1}.
$$
Then the Rankin-Selberg
$L$-function $L(s,\pi\otimes\pi')$ is a
Dirichlet series
\bea\label{4/RSL-}
L(s,\pi\otimes\pi')=
\sum_{n=1}^\infty \f{\l_{\pi\otimes\pi'}(n)}{n^s}
\eea
which is proved to be absolutely convergent for $\s>1$.
At the infinite place,
$$
L_\infty(s,\pi_\infty\otimes\pi'_\infty)
=\prod_{j=1}^{mm'}
\G_{\mathbb R}(s-\mu_{\pi\otimes\pi'}(j)).
$$
Moreover, at places $v$ where $\pi_v$ is unramified,
$L_v(s,\pi_v\otimes\pi'_v)$ has the following
explicit expression
\bea\label{4/Lp()=}
L_p(s,\pi_p\otimes\pi'_p)
=\prod_{j=1}^{m}\prod_{j'=1}^{m'}
\le(1-\f{\a_\pi(p,j)\a_{\pi'}(p,j')}{p^s}\ri)^{-1}
\eea
at $v=p$ a finite place, and at the
infinite place $v=\infty$,
\bea
L_\infty(s,\pi_\infty\otimes\pi'_\infty)
=\prod_{j=1}^{m}\prod_{j'=1}^{m'}
\G_{\mathbb R}(s-\mu_\pi(j)-\mu_{\pi'}(j')).
\eea
The complete $L$-function
$$
\Phi(s,\pi\otimes\pi')
=
L_\infty(s,\pi_\infty\otimes\pi'_\infty)
L(s,\pi\otimes\pi')
$$
satisfies a functional equation, and has properties
similar to those stated in the lemmas in \S4.1.
For simplicity, we do not list all these properties
of $L(s,\pi\otimes\pi')$ in detail, but just
point out some main differences between the Rankin-Selberg
$L$-function $L(s,\pi\otimes\pi')$ and the single $L$-function
$L(s,\pi)$:

\begin{itemize}
  \item The $\Phi(s,\pi\otimes\pi')$ has a meromorphic continuation
to $\mathbb C$;
  \item $\Phi(s,\pi\otimes\pi')$ is entire if
$\pi$ and $\pi'$ are not twisted equivalent, i.e.
$\pi'\not=\tilde\pi \otimes |\det|^{it}$ for any
$t\in \mathbb R$;
  \item if $\pi'=\tilde\pi \otimes |\det|^{it}$
for some $t\in \mathbb R$,
then
$L(s,\pi\otimes\pi')$ has only
a simple pole at $s=1+it$; in particular, the function
$L(s,\pi\otimes\tilde\pi)$ has only a simple pole at $s=1$.
\end{itemize}

For $\s>1$, the Euler product of
$L(s,\pi\otimes \pi')$ takes the form
\bea\label{RSLeuler/AUTO}
L(s,\pi\otimes \pi')
&=&\prod_{p}\prod_{j=1}^{m}\prod_{j'=1}^{m'}
\le(1-\f{\a_\pi(p,j)\a_{\pi'}(p,j')}{p^{s}}\ri)^{-1}.
\eea

The following result gives information for the Dirichlet
coefficients $\l_{\pi\otimes\tilde\pi}(n)$ for $L(s,\pi\otimes \tilde\pi)$;
for a proof of this, see Lemma A.1 in Rudnick and Sarnak \cite{RudSar}.

\noindent
{\bf Lemma 4.9.} (Rudnick-Sarnak \cite{RudSar}).
{\it Let $\pi$
an irreducible unitary cuspidal representation
for $GL_m(\mathbb A_\mathbb Q)$. Specifying $\pi'=\tilde{\pi}$ in (\ref{4/RSL-}), and write, for $\s>1$,
\bea\label{4/DEFLpipi}
L(s,\pi\otimes\tilde\pi)=\sum_{n=1}^\infty \f{\l_{\pi\otimes\tilde\pi}(n)}{n^s}.
\eea
Then
$$
\l_{\pi\otimes\tilde\pi}(n)\geq 0, \qquad \text{ for all } n\geq 1.
$$}

Other relations among the Dirichlet coefficients of
Rankin-Selberg $L$-functions $L(s,\pi\otimes\tilde\pi)$ will also be necessary.
We reserve the next section for this purpose.

\section{Coefficients of $L$-functions and Rankin-Selberg $L$-functions}

We need some general lemmas, which will be applied later to
Dirichlet coefficients of $L$-functions $L(s,\pi)$, or those for
Rankin-Selberg $L$-functions $L(s,\pi\otimes\tilde\pi)$.
The first general result is due to Brumley \cite{Bru},
and established by the theory of symmetric algebra.

\noindent
{\bf Lemma 4.10.} \text{(Brumley \cite{Bru}).}
{\it For $m$ complex numbers $\{\a_j\}_{j=1}^m$,
define the coefficients $b_n$ by
$$
\sum_{n=0}^\infty b_nX^n
=\prod_{j=1}^m \prod_{j'=1}^{m}
(1-\a_j\overline{\a_{j'}}X)^{-1}.
$$
If $\a_1\cdots\a_m=1$, then we have
$
b_m\geq 1.
$
In particular
for any irreducible unitary cuspidal representation of $GL_m({\mathbb A}_{\mathbb Q})$
and any prime $p$ such that $\pi_p$ is unramified, we have
$$\l_{\pi\otimes\tilde\pi}(p^m)\geq 1,$$
where $\l_{\pi\otimes\tilde\pi}(n)$ is defined by (\ref{4/DEFLpipi}).}

The second general lemma is due to L\"u \cite{Lu}. I am very
grateful for his kindness in allowing my reproduction of his proof below.

\noindent
{\bf Lemma 4.11.} (L\"u \cite{Lu}).
{\it For $m$ complex numbers $\{\a_j\}_{j=1}^m$,
define the coefficients $\ell_n$ by
\bea\label{4/cnDEF}
\sum_{n=0}^\infty \ell_nX^n
=\prod_{j=1}^m (1-\a_jX)^{-1}.
\eea
Also, for $n\geq 1$, define
\bea\label{4/dnDEF}
a_n=\a_1^n+\cdots+\a_m^n.
\eea
Then we have, for any $n\geq 1$,
\bea\label{4/ASSERTncn}
n \ell_n = a_1\ell_{n-1}+a_2\ell_{n-2}+\cdots
+a_{n-1}\ell_1+a_n.
\eea
}

\noindent{\it Proof.}
Differentiating (\ref{4/cnDEF}), we get
\bea\label{4/cnDEF+}
\sum_{n=1}^\infty n\ell_nX^{n-1}
&=&\sum_{i=1}^m \a_i(1-\a_iX)^{-1}\prod_{j=1}^m (1-\a_jX)^{-1}\nonumber\\
&=&\le(\sum_{i=1}^m \a_i(1-\a_iX)^{-1}\ri)\prod_{j=1}^m (1-\a_jX)^{-1}.
\eea
By expanding $(1-\a_iX)^{-1}$ and using the definition (\ref{4/dnDEF}),
the quantity within the last braces in (\ref{4/cnDEF+}) can be written as
\bna
\sum_{i=1}^m \a_i(1-\a_iX)^{-1}
&=&\sum_{i=1}^m \a_i \le(\sum_{u=0}^\infty \a_i^uX^u\ri)
\\
&=&\sum_{u=0}^\infty X^u\sum_{i=1}^m \a_i^{u+1}
\\
&=&\sum_{u=0}^\infty a_{u+1} X^u.
\ena
From this and (\ref{4/cnDEF}), one sees that the right-hand side in
(\ref{4/cnDEF+}) becomes
\bna
\le(\sum_{u=0}^\infty a_{u+1} X^{u}\ri)\le(\sum_{v=0}^\infty \ell_{v} X^{v}\ri)
&=&\sum_{n=0}^\infty \le(\sum_{\substack{u+v=n\\ u\geq 0, v\geq 0}} a_{u+1} \ell_v\ri)X^{n}
\\
&=&\sum_{n=0}^\infty \le(\sum_{\substack{u+v=n+1\\ u\geq 1, v\geq 0}} a_{u} \ell_v\ri)X^{n}
\\
&=&\sum_{n=1}^\infty \le(\sum_{\substack{u+v=n\\ u\geq 1, v\geq 0}} a_{u} \ell_v\ri)X^{n-1}.
\ena
Comparing this with the left-hand side of (\ref{4/cnDEF+}), we get, for all $n\geq 1$,
$$
n\ell_n=\sum_{\substack{u+v=n\\ u\geq 1, v\geq 0}} a_u \ell_v,
$$
which is exactly the assertion (\ref{4/ASSERTncn}) of the lemma.
\hfill $\square$

Applying the above two lemmas, we get the following consequence, which is
very important in establishing the the main result Theorem 4.15 of this chapter.

\noindent
{\bf Lemma 4.12.} {\it Let $m\geq 2$ be an integer and let $\pi$ be an irreducible unitary cuspidal
representation of $GL_m({\mathbb A}_{\mathbb Q})$.
For any prime $p$ such that $\pi_p$ is unramified, we have
$$
|\l_\pi(p^{m})|+|\l_\pi(p^{m-1})|+\cdots +|\l_\pi(p)|\geq \f{1}{m},
$$
where $\l_\pi(n)$ is as in (\ref{Lspi=Dir}) and (\ref{lpin=}).}

\noindent
{\it Proof.} The proof is divided into three steps, for a clear presentation.
The first two steps deal with the Rankin-Selberg $L$-function $L(s,\pi\otimes\tilde\pi)$
and the automorphic $L$-function $L(s,\pi)$, respectively, and the third is saved
for the final argument.

\medskip
\noindent
{\sc First step.} Let $\{\a_\pi(p,j)\}_{j=1}^m$ be the set of Satake
parameters for $\pi_p$; we may write $\a_\pi(p,j)=\a_j$ for simplicity.
Then (\ref{4/Lp()=}) becomes
\bea\label{4/Lp()=CON}
L_p(s,\pi_p\otimes\tilde\pi_p)
&=&\prod_{j=1}^{m}\prod_{j'=1}^{m}
\le(1-\f{\a_j\a_{j'}}{p^s}\ri)^{-1}\nonumber\\
&=&\prod_{\ell=1}^{M}
\le(1-\f{\b_\ell}{p^s}\ri)^{-1},
\eea
where we have written $M=m^2$ and
\bea\label{4/B=AA}
\{\b_\ell\}_{\ell=1}^M=\{\a_j\overline{\a_{j'}}\}_{1\leq j\leq m, 1\leq j'\leq m}.
\eea
Therefore, (\ref{4/Lp()=}) and (\ref{4/DEFLpipi}) give, for $\s>1$,
\bna\label{4/DEFLpipi+}
L_p(s,\pi_p\otimes \tilde\pi_p)
=\prod_{\ell=1}^{M}
\le(1-\f{\b_\ell}{p^{s}}\ri)^{-1}
=\sum_{n=0}^\infty \f{\l_{\pi\otimes\tilde\pi}(p^n)}{p^{ns}},
\ena
This is of the form (\ref{4/cnDEF}),
if we make the change of variables
$$
p^{-s}=X, \qquad \l_{\pi\otimes\tilde\pi}(p^n)=\ell_n.
$$
Thus, Lemma 4.11 gives, for all $n\geq 1$,
\bea\label{4/FORal}
n\l_{\pi\otimes\tilde\pi}(p^n)
&=&
a_{\pi\otimes\tilde\pi}(p)\l_{\pi\otimes\tilde\pi}(p^{n-1})
+a_{\pi\otimes\tilde\pi}(p^2)\l_{\pi\otimes\tilde\pi}(p^{n-2})
+\cdots \nonumber\\\noalign{\vskip 2mm}
&&
+a_{\pi\otimes\tilde\pi}(p^{n-1})\l_{\pi\otimes\tilde\pi}(p)
+a_{\pi\otimes\tilde\pi}(p^n),
\eea
where, because of (\ref{4/dnDEF}),
we have written
\bea\label{4/app=}
a_{\pi\otimes\tilde\pi}(p^n)=\b_1^{n}+\cdots+\b_M^{n}.
\eea
Now we write
\bea\label{4/REL/api}
a_{\pi}(p^n)=\a_1^{n}+\cdots+\a_m^{n};
\eea
then we have from (\ref{4/app=}) and (\ref{4/B=AA}) that
\bea\label{4/REL/apipi}
a_{\pi\otimes\tilde\pi}(p^n)=
\sum_{j=1}^m\sum_{j'=1}^m (\a_j\overline{\a_{j'}})^{n}
=|a_{\pi}(p^n)|^2\geq 0.
\eea
Inserting (\ref{4/REL/apipi}) into (\ref{4/FORal}), we get
\bea\label{4/FORal+}
n\l_{\pi\otimes\tilde\pi}(p^n)
&=&
|a_{\pi}(p)|^2\l_{\pi\otimes\tilde\pi}(p^{n-1})
+|a_{\pi}(p^2)|^2\l_{\pi\otimes\tilde\pi}(p^{n-2})
+\cdots \nonumber
\\\noalign{\vskip 2mm}
&&
+|a_{\pi}(p^{n-1})|^2\l_{\pi\otimes\tilde\pi}(p)
+|a_{\pi}(p^n)|^2,
\eea
which holds for all $n\geq 1$.
From it, we can deduce, by a simple induction on the integer $n$,
that $\l_{\pi\otimes\tilde\pi}(p^n)\geq 0$ for all $n\geq 1$.
This gives another proof of Lemma 4.9 due to Rudnick-Sarnak \cite{RudSar}, 
for those primes $p$ such that $\pi_p$ is unramified.

\medskip
\noindent
{\sc Second step.} Let $\{\a_j\}_{j=1}^m$ be as before.
Then (\ref{Lpspip}) becomes
\bna
L_p(s, \pi_p )=\prod_{j=1}^m
\le(1-\f{\a_j}{p^{s}}\ri)^{-1}.
\ena
Therefore, (\ref{Lspi=Dir}) gives, for $\s>1$,
$$
L_p(s,\pi_p)
=\prod_{j=1}^{m}
\le(1-\f{\a_j}{p^{s}}\ri)^{-1}
=\sum_{n=0}^\infty \f{\l_{\pi}(p^n)}{p^{ns}},
$$
with $\l_\pi(n)$ defined by (\ref{lpin=}).
On changing variables
$$
p^{-s}=X, \qquad \l_{\pi}(p^n)=\ell_n,
$$
the above is again of the form (\ref{4/cnDEF}), and
Lemma 4.11 gives, for all $n\geq 1$,
\bea\label{4/FORal++}
n\l_{\pi}(p^n)
&=&
a_{\pi}(p)\l_{\pi}(p^{n-1})
+a_{\pi}(p^2)\l_{\pi}(p^{n-2})
+\cdots \nonumber
\\\noalign{\vskip 2mm}
&&
+a_{\pi}(p^{n-1})\l_{\pi}(p)
+a_{\pi}(p^n),
\eea
where $a_\pi(p^n)$ is as in (\ref{4/REL/api}), as suggested by (\ref{4/dnDEF}).

\medskip
\noindent
{\sc Third step.} Taking $n=m$ in (\ref{4/FORal+}), we have
\bea\label{4/FORal+++}
m\l_{\pi\otimes\tilde\pi}(p^m)
&=&
|a_{\pi}(p)|^2\l_{\pi\otimes\tilde\pi}(p^{m-1})
+|a_{\pi}(p^2)|^2\l_{\pi\otimes\tilde\pi}(p^{m-2})
+\cdots \nonumber
\\\noalign{\vskip 2mm}
&&
+|a_{\pi}(p^{m-1})|^2\l_{\pi\otimes\tilde\pi}(p)
+|a_{\pi}(p^m)|^2.
\eea
Brumley's lemma (Lemma 4.10) now asserts that
$\l_{\pi\otimes\tilde\pi}(p^m)\geq 1,$
and therefore, the left-hand side in (\ref{4/FORal+++})
is bounded from below by
$m\l_{\pi\otimes\tilde\pi}(p^m)\geq m.$
As we have seen that  $\l_{\pi\otimes\tilde\pi}(p^n)\geq 0$ for all $n\geq 1$,
each term on the right-hand side of (\ref{4/FORal+++})
is non-negative. These observations will be used later.

Before going further, we make a claim:

\noindent
{\sc Claim C.} {\it There exists a positive integer $j$ with $1\leq j\leq m$ such that
\bea\label{4/CLAIM}
|a_{\pi}(p^j)|^2\geq 1.
\eea}

We suppose that Claim C is not true, and establish a contradiction.
Since Claim C is not true, we must have
\bea\label{4/INVc}
|a_{\pi}(p^n)|^2<1
\eea
for all $1\leq n\leq m$.
It thus follows from (\ref{4/FORal+++}) and (\ref{4/INVc}) with $n=1$ that
\bea\label{4/l=a<1}
\l_{\pi\otimes\tilde\pi}(p)=|a_{\pi}(p)|^2<1.
\eea
We may also apply (\ref{4/FORal+++}) and (\ref{4/INVc}) with $n=2$,
and we get from (\ref{4/l=a<1}) that
\begin{align*}
2\l_{\pi\otimes\tilde\pi}(p^2)
& =
|a_{\pi}(p)|^2\l_{\pi\otimes\tilde\pi}(p)
+|a_{\pi}(p^2)|^2
\\\noalign{\vskip 1mm}
& <1+1=2,
\end{align*}
that is $\l_{\pi\otimes\tilde\pi}(p^2)<1.$
By induction on the $n$ in (\ref{4/FORal+}),
we can prove that
$\l_{\pi\otimes\tilde\pi}(p^n)<1$ for all $1\leq n\leq m$.
In particular,
$
\l_{\pi\otimes\tilde\pi}(p^m)<1.
$
This contradicts Brumley's lemma (Lemma 4.10), which asserts that
$$\l_{\pi\otimes\tilde\pi}(p^m)\geq 1$$
in the present situation.
Hence, Claim C is proved.

\medskip

If $n_0$ is one of the integers in Claim C so that (\ref{4/CLAIM}) holds,
then, by (\ref{4/REL/apipi}), one has $|a_{\pi}(p^{n_0})|^2\geq 1$,
and consequently,
\bea\label{4/apn0>1}
|a_{\pi}(p^{n_0})|\geq 1.
\eea
Now let $n_0$ with $1\leq n_0\leq m$ be the smallest integer such that
(\ref{4/apn0>1}) holds. It follows that
$$
|a_\pi(p^j)|<1, \qquad \text{ for all } 1\leq j<n_0.
$$
Applying (\ref{4/FORal++}) with $n=n_0$, we deduce that
\bna
n_0|\l_\pi(p^{n_0})|
&=&|a_{\pi}(p)\l_{\pi}(p^{n_0-1})
+\cdots \nonumber
+a_{\pi}(p^{n_0-1})\l_{\pi}(p)
+a_{\pi}(p^{n_0})|
\\\noalign{\vskip 1mm}
&>& -|\l_{\pi}(p^{n_0-1})|
-\cdots -|\l_{\pi}(p)|
+1.
\ena
This implies that
\bna
m\{|\l_\pi(p^m)|+\cdots +|\l_\pi(p)|\}
&\geq& m\{|\l_\pi(p^{n_0})|+\cdots +|\l_\pi(p)|\}
\\\noalign{\vskip 1mm}
&\geq& n_0 |\l_\pi(p^{n_0})|+\cdots +|\l_\pi(p)|
\\\noalign{\vskip 1mm}
&>& 1,
\ena
and the lemma follows.
\hfill $\square$

\section{Sign changes of $\l_\pi(n)$}

Working similarly as in Knopp, Kohnen, and Pribitkin \cite{KnoKohPri}, we can prove

\noindent
{\bf Theorem 4.13.} {\it Let $m\geq 2$ be an integer and let $\pi$ be an irreducible unitary cuspidal representation
for $GL_m(\mathbb A_\mathbb Q)$
such that $\l_\pi(n)$ are real for all $n\geq 1$. Then the sequence
$
\{\l_f(n)\}_{n=1}^\infty
$
has infinitely many sign changes, i.e. there are
infinitely many $n$ such that $\l_f(n)>0$, and there are infinitely
many $n$ such that $\l_f(n)<0$.
}

\noindent
{\bf Corollary 4.14.} {\it Let $m\geq 2$ be an integer and
let $\pi$ be an irreducible unitary cuspidal representation
for $GL_m(\mathbb A_\mathbb Q)$
such that it is self-contragredient. Then $\l_\pi(n)$ are real for all $n\geq 1$, and
the sequence
$
\{\l_f(n)\}_{n=1}^\infty
$
has infinitely many sign changes, i.e. there are
infinitely many $n$ such that $\l_f(n)>0$, and there are infinitely
many $n$ such that $\l_f(n)<0$.
}

\section{A Linnik-type problem for automorphic $L$-functions}

Now we state our main result in this chapter.

\noindent
{\bf Theorem 4.15.}
{\it Let $m\geq 2$ be an integer and
let $\pi$ be an irreducible unitary cuspidal representation for $GL_m(\mathbb A_\mathbb Q)$.
If $\l_\pi(n)$ are real for all $n\geq 1$, then there is some $n$ satisfying
\bea\label{bound/AUTO}
n\ll Q_\pi^{m/2+\ve}
\eea
such that $\l_\pi(n)<0$. The constant implied in (\ref{bound/AUTO})
depends only on $m$ and $\ve$. In particular, the result is true
for any self-contragredient representation $\pi$. }

\noindent
{\it Proof. } Still, let us start with the sum
$$
S(x):=\sum_{n\leq x} \l_\pi(n)\bigg(\log \f{x}{n}\bigg)^\ell,
$$
assuming that
\bea\label{assum/AUTO}
\l_\pi(n)\geq 0, \quad \text{\rm for} \ n\leq x.
\eea
Here $\ell$ is a positive integer that will be decided later.
The desired result will follow from upper and lower bound estimates
for $S(x)$.

To get an upper bound for $S(x)$, we apply Perron's formula (\ref{1/perSIM})
to the Dirichlet series (\ref{Lspi=Dir}), getting
$$
S(x)=\f{1}{2\pi i}\int_{2-i\infty}^{2+i\infty}L(s,\pi)\f{x^s}{s^{\ell+1}}{\rm d} s.
$$
Moving the contour to the vertical line $\s=\ve$ with $\ve$ being an arbitrarily small
positive constant, and applying Harcos'
convexity bound (\ref{conv/AUTO}) for $L(s,\pi)$, we obtain
\bna
S(x)&=&\f{1}{2\pi i}\int_{\ve-i\infty}^{\ve+i\infty}L(s,\pi)\f{x^s}{s^{\ell+1}}{\rm d} s \\
&\ll_\ve & \int_{-\infty}^{\infty} Q_\pi(t)^{1/2+\ve}\f{x^{\ve}}{(|t|+\ve)^{\ell+1}}{\rm d} t.
\ena
The analytic conductor $Q_\pi(t)$ is bounded by $|t|^m$ in the $t$-aspect.
Thus, if we take $\ell=m$, then the above estimate becomes
\bea\label{5/Supp}
S(x)
&\ll_{\ell,m,\ve}& Q_\pi^{1/2+\ve}x^\ve \int_{-\infty}^{\infty}
\f{(|t|+1)^{m/2}}{(|t|+\ve)^{\ell+1}}{\rm d} t\nonumber\\
&\ll_{m,\ve}& Q_\pi^{1/2+\ve}x^\ve.
\eea
This is the desired upper bound for $S(x)$.

To get a lower bound for $S(x)$, we apply Lemma 4.12, from which we have
\bna
S(x)
&\geq& (\log 2)^\ell \sum_{\substack{n\leq x/2\\ (n,N_\pi)=1}} \l_\pi(n)
\\
&\geq& (\log 2)^\ell \sum_{\substack{p\leq (x/2)^{1/m}\\ p\nmid N_\pi}}
\{\l_\pi(p^m)+\l_\pi(p^{m-1})+\cdots+\l_\pi(p)\}\\
&\gg_{\ell, m}& \sum_{\substack{p\leq (x/2)^{1/m}\\ p\nmid N_\pi}} 1
\gg_{\ell, m} \f{(x/2)^{1/m}}{\log (x/2)}-\log N_\pi.
\ena
Without loss of generality, we may suppose
\bea\label{4/DXBIG}
x\geq C\log^{m+1}Q_\pi;
\eea
where $C$ is a large absolute constant.
This requirement is very mild in view of the assertion of the theorem.
\hfill $\square$

\chapter{The prime number theorem for automorphic $L$-functions}

\section{The automorphic prime number theorem}

\noindent
To each irreducible unitary cuspidal representation
$\pi=\otimes \pi_p$ of $GL_m({\mathbb A}_{\mathbb Q})$, one can attach a global $L$-function
$L(s,\pi)$ as in \S 4.1.
Then, one can link $L(s,\pi)$ with primes by taking logarithmic differentiation
in (\ref{Lpspip}), so that for $\s>1$,
\bea\label{5/L'/L}
\frac{{\rm d}}{{\rm d}s} \log L(s,\pi)
=
-\sum_{n=1}^\infty
\frac{\L(n)a_\pi(n)}{n^s},
\eea
where $\L(n)$ is the von Mangoldt function, and
\bea\label{apipk/AUTO}
a_\pi(p^k)=
\sum_{j=1}^m
\alpha_\pi(p,j)^k.
\eea
It will be important later that this is the same as that defined in
(\ref{4/REL/api}).
The prime number theorem for $L(s,\pi)$ concerns the asymptotic
behavior of the counting function
\bna
\psi(x,\pi)=\sum_{n\leq x}\L(n)a_\pi(n),
\ena
and a special case of it
asserts that, if $\pi$ is an irreducible unitary cuspidal
representation of $GL_m({\mathbb A}_{\mathbb Q})$ with $m\geq 2$, then
\bea\label{psixpi/AUTO}
\psi(x,\pi)\ll \sqrt{Q_\pi} \cdot x \cdot \exp\bigg(-\f{c}{2m^4}\sqrt{\log x}\bigg)
\eea
for some absolute positive constant $c$, where the implied constant is absolute.
In Iwaniec and Kowalski \cite{IwaKow}, Theorem 5.13, a prime number theorem is proved for general
$L$-functions satisfying necessary axioms, from which (\ref{psixpi/AUTO})
follows as a consequence.

In this chapter, we investigate the influence of GRH on $\psi(x,\pi)$.
It is known that, under GRH, (\ref{psixpi/AUTO}) can be improved to
\bea\label{undGRH/AUTO}
\psi(x,\pi)\ll x^{1/2}\log^2 (Q_{\pi}x),
\eea
but better results are desirable. For a proof of (\ref{undGRH/AUTO}),
see e.g. \S5.3.
In this direction, we establish the following results.

\noindent
{\bf Theorem 5.1.} {\it Let $m\geq 2$ be an integer and let $\pi$ be an irreducible unitary cuspidal
representation of $GL_m({\mathbb A}_{\mathbb Q})$. Assume GRH for $L(s,\pi)$. Then
we have
\bea\label{undGRH/AUTO+}
\psi(x,\pi)\ll x^{1/2}\log^2(Q_{\pi}\log x)
\eea
for $x\geq 2$, except on a set $E$ of finite logarithmic measure, i.e.
\bna
\int_{E} \f{{\rm d}x}{x}< \infty.
\ena
The constant implied in the
$\ll$-symbol depends
at most on $m$.}

Theorem 5.1 tells that, except on a set of finite logarithm measure $E$,
(\ref{undGRH/AUTO}) can be improved to (\ref{undGRH/AUTO+}). The next two theorems
say that $\psi(x,\pi)$ behaves somehow like $x^{1/2}$.

\noindent
{\bf Theorem 5.2.} {\it Let $m\geq 2$ be an integer and let $\pi$ be an irreducible unitary cuspidal
representation of $GL_m({\mathbb A}_{\mathbb Q})$. Assume GRH for $L(s,\pi)$.
Then
\bna
\int^X_2 |\psi(x,\pi)|^2\f{{\rm d}x}{x}\ll X\log^2Q_{\pi}.
\ena
The constant implied in the
$\ll$-symbol depends
at most on $m$.}

\goodbreak

Gallagher \cite{Gal1} was the first to establish a result
like Theorem 5.1, in the classical case $m=1$ for the Riemann
zeta-function. He proved that, under the Riemann Hypothesis for
the classical zeta-function,
\begin{align*}
\psi(x)
& :=\sum_{n\leq x}\L(n)
=x+O\big(x^{1/2}(\log \log x)^2\big)
\end{align*}
for $x\geq 2$, except on a set of finite logarithmic measure,
and hence made improvement on the classical estimate error term
$O(x^{1/2}\log^2 x)$ of von Koch \cite{Koc}. In the same paper, Gallagher \cite{Gal1}
also gave short proofs of Cram\'er's conditional estimate (see \cite{Cra1} and \cite{Cra2})
\bna
\int_2^X (\psi(x)-x)^2\f{{\rm d}x}{x}\ll X.
\ena
Gallagher's
proofs of the above results make crucial use of his lemma in \cite{Gal2},
which is now named after him.

Our Theorems 5.1-5.2 generalize the above classical results to
the prime counting function $\psi(s,\pi)$ attached to
irreducible unitary cuspidal
representations $\pi$ of $GL_m({\mathbb A}_{\mathbb Q})$ with $m\geq 2$.
Our proofs combine the approach of Gallagher with recent results
of Liu and Ye (\cite{LiuYe1}, \cite{LiuYe3})
on the distribution of zeros of Rankin-Selberg automorphic $L$-functions.

The above Theorem 5.2 states that, under GRH, $|\psi(x,\pi)|$ is of size $x^{1/2}\log Q_\pi$
on average. This can be compared with the next theorem, which gives
the unconditional Omega result that $|\psi(x,\pi)|$ should not be of order lower
than $x^{1/2-\ve}$.

\noindent
{\bf Theorem 5.3.} {\it Let $m\geq 2$ be an integer and let $\pi$ be an irreducible unitary cuspidal
representation of $GL_m({\mathbb A}_{\mathbb Q})$, and $\ve>0$ arbitrary.
Unconditionally,
\bna
\psi(x,\pi)=\Omega(x^{1/2-\ve}),
\ena
where the constant implied in the
$\O$-symbol depends
at most on $m$ and $\ve$. More precisely,
there exists an increasing sequence $\{x_n\}_{n=1}^\infty$ tending to infinity such that
\bea\label{5/LIMd}
\lim_{n\to\infty}
\f{|\psi(x_n,\pi)|}{ x_n^{1/2-\ve}}>0.
\eea
}

Note that the sequence $\{x_n\}_{n=1}^\infty$ and the limit in (\ref{5/LIMd})
may depend on $\pi$ and $\ve$.
This result generalizes that for the Riemann zeta-function. It is possible to get
better Omega results like those in Chapter V of Ingham \cite{Ing}.
We remark that, unlike the classical case, in Theorems 5.1-5.3 we do not have
the main term $x$. This is because $L(s,\pi)$
is entire when $m\geq 2$, while $\zeta(s)$ has a simple pole at
$s=1$ with residue $1$.

\section{Preliminaries}

We need some preliminaries to establish the main results.

\noindent
{\bf Lemma 5.4.} {\it Let $\pi$ be an irreducible unitary cuspidal
representation of $GL_m({\mathbb A}_{\mathbb Q})$ with $m\geq 2$. Then
\bna
\f{{\rm d}}{{\rm d}s}\log L(s,\pi)
&=& C+\sum_{\rho}\le(\f{1}{s-\rho}+\f{1}{\rho}\ri)
+\sum_{j=1}^m \f{1}{s+\mu_{\pi}(j)}\\
&& +\sum_{j=1}^m \sum_{n=1}^{\infty}\le(\f{1}{2n+s+\mu_{\pi}(j)}-\f{1}{2n}\ri),
\ena
where $C$ is a constant depending on $\pi$. The set of all trivial zeros of $L(s,\pi)$ is
$$
\{\mu: \mu=-2n-\mu_{\pi}(j), \quad n=0,1,2,\dots; \ j=1,\ldots,m\}.
$$}

\noindent
{\it Proof.} Since $\Phi(s,\pi)$ is of order one (Lemma 4.4),
we have (see e.g. Davenport \cite{Dav}, Chapter 11)
$$
\Phi(s,\pi)=e^{A+Bs}\prod_{\rho}\le(1-\f{s}{\rho}\ri)e^{s/\rho},
$$
where $A,B$ are constants depending on $\pi.$ Taking logarithmic
derivative, we get
\bea\label{dds/AUTO}
\f{{\rm d}}{{\rm d}s}\log \Phi(s,\pi)
=B+\sum_{\rho}\le(\f{1}{s-\rho}+\f{1}{\rho}\ri),
\eea
where we set $\log 1=0.$ By the definition of
$\Phi(s,\pi)$,
\bea\label{dds/AUTO+}
\f{{\rm d}}{{\rm d}s}\log \Phi(s,\pi)
=\f{{\rm d}}{{\rm d}s}\log L_{\infty}(s,\pi_{\infty})
+\f{{\rm d}}{{\rm d}s}\log L(s,\pi).
\eea
Applying
$$
\f{{\rm d}}{{\rm d}s}\log \Gamma(s)=-\f{1}{s}-\g
-\sum_{n=1}^{\infty}\le(\f{1}{n+s}-\f{1}{n}\ri),
$$
where $\g$ is Euler's constant, we have
\bna
\f{{\rm d}}{{\rm d}s}\log L_{\infty}(s,\pi_{\infty})
&=&\sum_{j=1}^m \f{{\rm d}}{{\rm d}s}\log \pi^{-(s+\mu_{\pi}(j))/2}
+\sum_{j=1}^m \f{{\rm d}}{{\rm d}s}\log \Gamma\le(\f{s+\mu_{\pi}(j)}{2}\ri)\\
&=& -\f{m}{2}(\log \pi +\g)
-\sum_{j=1}^m \f{1}{s+\mu_{\pi}(j)}\\
&&-\sum_{j=1}^m \sum_{n=1}^{\infty}\le(\f{1}{2n+s+\mu_{\pi}(j)}-\f{1}{2n}\ri).
\ena
Inserting this and (\ref{dds/AUTO}) into (\ref{dds/AUTO+}), we get the lemma. \hfill $\square$

\medskip

Consider the poles of
\bea\label{L1-s=G}
L(1-s,\tilde\pi_\infty)
=\pi^{-ms/2}\prod_{j=1}^{m}
\Gamma\left(\frac{1-s+\mu_{\tilde\pi}(j)}{2}\right).
\eea
These poles are easily to be seen as
$$
\{P_{n,j}=2n+1+\mu_{\tilde\pi}(j):\ n=0,1,2,\ldots, \quad j=1,\ldots,m\}.
$$
As in \cite{LiuYe1}, we let ${\mathbb C}(m)$ denote the complex plane with the discs
$$
|s-P_{n,j}|<\f{1}{8m},
\quad n=0,1,2,\ldots, \quad j=1,\ldots,m
$$
excluded. Thus, for any $s\in {\mathbb C}(m)$, the quantity
$$
\frac{1-s+\mu_{\tilde\pi}(j)}{2}
$$
is away from all poles of  $\Gamma(s)$ by at least $1/(16m)$.
Now we give a remark about the structure of ${\mathbb C}(m)$. 
For $j=1,\ldots,m$, denote by $\beta(j)$ the fractional
part of $\Re e \mu_{\tilde\pi}(j)$. In addition we let
$\beta(0)=0$ and $\beta(m+1)=1$.
Then all $\beta(j)\in [0,1]$, and hence there exist
$\beta(j_1),\beta(j_2)$ such that
$\beta(j_2)-\beta(j_1)\geq 1/(3m)$ and there is no
$\beta(j)$ lying between $\beta(j_1)$ and $\beta(j_2)$. It
follows that the strip
$$
S_0=\{s:\beta(j_1)+1/(8m)\leq \Re s\leq\beta(j_2)-1/(8m)\}
$$
is contained in ${\mathbb C}(m).$ Consequently, for all $n=0,1,2,\ldots,$
the strips
\bea\label{DEF/Sn}
S_n=\bigg\{s: -n+\beta(j_1)+\f{1}{8m}\leq \Re s\leq -n+\beta(j_2)-\f{1}{8m}\bigg\}
\eea
are subsets of ${\mathbb C}(m).$ In \cite{LiuYe1}, \S 4, Liu and Ye studied distribution of zeros of the
Rankin-Selberg $L$-function $L(s,\pi\otimes\pi')$, where $\pi$ and $\pi'$ are irreducible unitary cuspidal
representations of $GL_m({\mathbb A}_{\mathbb Q})$ and $GL_{m'}({\mathbb A}_{\mathbb Q})$,
respectively. This structure of ${\mathbb C}(m)$ is a special case of the
${\mathbb C}(m,m')$ in \cite{LiuYe1}, \S 4.

The following Lemma 5.5(i) and (ii)
are Lemma 4.3(d) and Lemma 4.4 of \cite{LiuYe1}, respectively.

\noindent
{\bf Lemma 5.5.} {\it Let $\pi$ be an irreducible unitary cuspidal
representation of $GL_m({\mathbb A}_{\mathbb Q})$ with $m\geq 2$. Then

\vskip -2mm

{\rm (i)} For $|T|>2,$ there exists $\tau$
with $T\leq \tau\leq T+1$ such that when $-2\leq \s\leq 2$,
$$
\f{{\rm d}}{{\rm d}s}\log L(\s \pm i\tau,\pi)\ll
\log^2 (Q_{\pi}|\tau|).
$$

\vskip -2mm

{\rm (ii)} If $s$ is in some strip $S_n$ as in $(\ref{DEF/Sn})$ with
$n\leq -2,$ then
$$
\f{{\rm d}}{{\rm d}s}\log L(s,\pi) \ll 1.
$$}

\section{An explicit formula}

The explicit formula given in Theorem 5.6 below
is unconditional; it requires
neither GRH nor GRC.

\noindent
{\bf Theorem 5.6.} {\it Let $m\geq 2$ be an integer and let $\pi$ be an irreducible unitary cuspidal
representation of $GL_m({\mathbb A}_{\mathbb Q})$. Then, for
$x\geq 2$ and $T\geq 2$,
\bna
\psi(x,\pi)
&=&-\sum_{|\g|\leq T}\f{x^\rho}{\rho}
+O\le\{\min \le(\f{x}{T^{1/4}}, \f{x^{1+\th}}{T^{1/2}}\ri)\log (Q_\pi x)\ri\} \\
&& +O(x^\th \log x)
+O\le(\f{x\log^2 (Q_{\pi}x)}{T^{1/2}}\ri) + O\le(\f{\log T}{x}\ri),
\ena
where $\th$ is as in Lemma 4.8.}

We will establish Theorem 5.6 at the end of this section.
Explicit formulas of different forms were
established by Moreno (\cite{Mor1},  \cite{Mor2}); under GRC, explicit formulas
for general $L$-functions were proved in (5.53) of Iwaniec and
Kowalski \cite{IwaKow}.

Using Theorem 5.6, we will by the way give a proof of the prime number theorem
for $\psi(x,\pi)$ under GRH.

\noindent
{\bf Corollary 5.7.} {\it Let $m\geq 2$ be an integer and let $\pi$ be an irreducible unitary cuspidal
representation of $GL_m({\mathbb A}_{\mathbb Q})$, and assume
GRH for $L(s,\pi)$. Then (\ref{undGRH/AUTO}) holds.}

\noindent
{\it Proof.} Theorem 5.6 with $T=x^8$ gives
$$
\psi(x,\pi)=-\sum_{|\g|\leq T}\f{x^\rho}{\rho}
+O(x^\th \log x)+O\bigg(\f{\log^2(Q_{\pi}x)}{x}\bigg).
$$
By Lemma 4.8, the error term is acceptable in (\ref{undGRH/AUTO}).
Under GRH for $L(s,\pi)$, we have $\rho=1/2+i\g$ in the formula
above, and therefore, by partial summation and Lemma 4.6,
\bna
\sum_{|\g|\leq T}\f{x^\rho}{\rho}
&\ll& x^{1/2} \bigg(\sum_{|\g|\leq 1} 1 + \sum_{1\leq |\g|\leq T}\f{1}{|\gamma|}\bigg)
\\
&\ll& x^{1/2}\bigg\{\log Q_\pi+\int_1^T \f{1}{t}{\rm d}N(t,\pi)\bigg\}
\\\noalign{\vskip 2mm}
&\ll& x^{1/2}\log^2 (Q_{\pi}T).
\ena
This proves (\ref{undGRH/AUTO}).
\hfill $\square$

\medskip

The following form of Perron's formula will be needed in
the proof of Theorem 5.6.

\noindent
{\bf Lemma 5.8.} (Perron's formula).
{\it Under the assumption of Lemma 3.11,
we have, for $b>\s_a, x\geq 2$, $T\geq 2$,
\bea\label{perron/NOR}
\sum_{n\leq x}a_n
&=&\f1{2\pi i}\int_{b-iT}^{b+iT}F(s)\f{x^{s}}{s}{\rm d} s
+O\bigg(\sum_{|n-x|\leq x/\sqrt{T}} |a_n|+ \f{x^{b}B(b)}{\sqrt{T}}\bigg).
\eea}

A key feature of Lemma 5.8 is that individual upper bound for $a_n$ does not appear
on the right-hand side, and this makes Theorem 5.6, and hence Theorems 5.1-5.3,
independent of GRC.  Perron's formula of this nature was successfully applied in classical cases
where $a_n$ is not bounded by $1$. The specific form of Perron's formula in Lemma 5.8, though not optimal,
must have been known to the expert for some time. It follows
from Tenenbaum \cite{Ten}, Theorem II.2.2,
for example. See also Liu and Ye \cite{LiuYe3} for a proof and some applications to
automorphic $L$-functions.

\noindent
{\bf Lemma 5.9.} (Iwaniec-Kowalski \cite{IwaKow}).
{\it Let $\pi$ be an irreducible unitary cuspidal
representation of $GL_m({\mathbb A}_{\mathbb Q})$ with $m\geq 2$. Then
$$
\sum_{n\leq x}|\Lambda(n)a_\pi(n)|^2\ll m^2 x \log^2(Q_\pi x),
$$
where the implied constant is absolute. }

This is (5.48) in \cite{IwaKow}, and proved in the lower part on page 110
of \cite{IwaKow}.

\noindent
{\it Proof of Theorem 5.6.} In view of (\ref{Lspi=Dir}) and Lemma 4.1,
we can apply Lemma 5.8 with $\s_a=1, b=1+1/\log x$, and
$$
F(s)=\f{{\rm d}}{{\rm d}s}\log L(s,\pi)=-\sum_{n=1}^\infty
\frac{\L(n)a_\pi(n)}{n^s},
$$
that is $a_n=-\L(n)a_\pi(n)$.

To estimate the first $O$-term in (\ref{perron/NOR}), we let $0< y\leq x$, and consider
\bna
\sum_{x<n\leq x+y} |\L(n)a_\pi(n)|
&\ll&
\biggl\{\sum_{n\leq 2x}|\L(n)a_\pi(n)|^2\biggr\}^{1/2}
\biggl\{\sum_{x<n\leq x+y}1\biggr\}^{1/2}
\\
&\ll& \sqrt{x(y+1)}\log (Q_\pi x).
\ena
On the other hand, by (\ref{apipk/AUTO}) and the bound toward GRC in Lemma 4.8, for $n=p^k$,
$$
|a_\pi(n)|=|a_\pi(p^k)|\leq \sum_{j=1}^m |\a_\pi (p,j)|^k \leq mp^{k\th}\leq mn^\th.
$$
Therefore, trivially,
$$
\sum_{x<n\leq x+y} |\L(n)a_\pi(n)|\ll x^{\th}(y+1)\log x.
$$
Hence,
\bea\label{sec/1}
&&\sum_{|n-x|\leq x/\sqrt{T}} |\L(n)a_\pi(n)| \nonumber\\
&&\ll \min \le\{\sqrt{x\le(\f{x}{T^{1/2}}+1\ri)}\log (Q_\pi x), \
x^\th \le(\f{x}{T^{1/2}}+1\ri)\log x\ri\}\nonumber\\
&&\ll \min \le(\f{x}{T^{1/4}}, \f{x^{1+\th}}{T^{1/2}}\ri)\log (Q_\pi x) + x^\th \log x.
\eea
In the last step, we have considered the two cases $T\leq x^2$ and $T> x^2$
separately. The other $O$-term in (\ref{perron/NOR}) depends on
$B(\s)$. For $\s>1$, Cauchy's inequality gives
\bna
B(\s)=\sum_{n=1}^\infty \f{|\L(n)a_\pi(n)|}{n^\s}
\ll
\biggl\{\sum_{n=1}^\infty \f{|\L(n)a_\pi(n)|^2}{n^\s}\biggr\}^{1/2}
\biggl\{\sum_{n=1}^\infty \f{1}{n^\s}\biggr\}^{1/2}.
\ena
By Lemma 5.9,
\bna
\sum_{n=1}^\infty \f{|\L(n)a_\pi(n)|^2}{n^\s}
&=&
\int_1^\infty \f 1{u^\s} {\rm d} \biggl\{\sum_{n\leq u}|\L(n)a_\pi(n)|^2\biggr\}
\nonumber\\
&\ll& \log^2Q_\pi
+\int_1^\infty \f{\log^2 (Q_\pi u)}{u^\s}{\rm d} u \nonumber\\
&\ll& \f{\log^2Q_\pi}{\s-1}+\f{1}{(\s-1)^3}.
\ena
Similarly but more easily, we have
$$
\sum_{n=1}^\infty \f{1}{n^\s}\ll \f{1}{\s-1},
$$
and consequently,
\bna
B(\s)\ll\f{\log Q_\pi}{\s-1}+\f{1}{(\s-1)^2}.
\ena
Therefore, the second $O$-term in (\ref{perron/NOR}) is
\bea\label{sec/2}
\ll \f{x(\log x)(\log Q_\pi x)}{\sqrt{T}}.
\eea
Inserting (\ref{sec/2}) and (\ref{sec/1}) into (\ref{perron/NOR}), we get
\bea\label{perron/NOR+}
\sum_{n\leq x}\L(n)a_\pi(n)
&=&\f{1}{2\pi i}
\int_{b-iT}^{b+iT}\le\{-\f{L'}L(s,\pi)\ri\}\f{x^s}{s}{\rm d} s\nonumber\\
&&+O\le\{\min \le(\f{x}{T^{1/4}}, \f{x^{1+\th}}{T^{1/2}}\ri)\log (Q_\pi x)\ri\} \nonumber\\
&&+O\bigg(\f{x(\log x)(\log Q_\pi x)}{\sqrt{T}}\bigg)
+O(x^\th \log x).
\eea

Next, we shall shift the contour of integration to the left.
Choose $a$ with $-2<a< -1$ such that the vertical line $\s=a$ is contained
in the strip $S_{-2}\subset {\mathbb C}(m)$; this is guaranteed by
the structure of ${\mathbb C}(m)$. Without loss of generality, let
$T>0$ be a large number such that $T$ and $-T$ can be taken as the
$\tau$ in Lemma 5.5(i). Now we consider the contour $C_1 \cup C_2 \cup C_3$
with
\bna
C_1: &
\quad b\geq \s\geq a, \ \ t=-T;
\\
C_2: & \quad \s=a, \ \ -T\leq t\leq T;
\\
C_3: & \quad a\leq \s\leq b, \ \ t=T.
\ena
By Lemma 5.4, certain nontrivial zeros $\rho=\b+i\g$ and trivial zeros $\mu=\l+i\nu$
of $L(s,\pi)$,
as well as $s=0$ are passed by the
shifting of the contour. Computing the residues by Lemma 5.4, we have
\bea\label{int/CON}
\hskip -8mm
\f{1}{2\pi i}
\int_{b-iT}^{b+iT}\le\{-\f{L'}L(s,\pi)\ri\}\f{x^s}{s}{\rm d} s
&=&\f{1}{2\pi i}
\le(\int_{C_1}+\int_{C_2}+\int_{C_3}\ri)\nonumber\\
&&
-\sum_{|\g|\leq T}\f{x^\rho}{\rho}
-\sum_{\substack{a<-\l<b\\ |\nu|\leq T}} \f{x^{-\mu}}{-\mu}-\f{L'}L(0,\pi).
\eea

The integral on $C_1$ can be estimated by Lemma 5.5(i) as
\begin{align*}
\f{1}{2\pi i} \int_{C_1} \le\{-\f{L'}L(s,\pi)\ri\}\f{x^s}{s} {\rm d}s
& \ll \int_{a}^{b} \log^2 (Q_{\pi}T) \f{x^\s}{T} {\rm d}\sigma
\\
& \ll \f{x\log^2 (Q_{\pi}T)}{T},
\end{align*}
and the same upper bound also holds for the integral on $C_3$. By
Lemma 5.5(ii), then
\begin{align*}
\f{1}{2\pi i} \int_{C_2} \le\{-\f{L'}L(s,\pi)\ri\}\f{x^s}{s} {\rm d}s
& \ll \int_{-T}^{T}\f{x^{a}}{|t|+1}{\rm d} t
\\
& \ll x^a \log T.
\end{align*}
To bound the contribution from the trivial zeros $\mu=\l+i\nu$, we apply
Lemma 4.8, so that
$$
\sum_{\substack{a<-\l<b\\ |\nu|\leq T}}\f{x^{-\mu}}{\mu}\ll x^\th,
$$
where we have used the fact that there are finite number of
trivial zeros $\mu=\l+i\nu$ with $a<-\l<b, |\nu|\leq T.$
Therefore, (\ref{int/CON}) becomes
\bna
&& \f{1}{2\pi i}
\int_{b-iT}^{b+iT}\le\{-\f{L'}L(s,\pi)\ri\}\f{x^s}{s}{\rm d} s \\
&& =-\sum_{|\g|\leq T}\f{x^\rho}{\rho}
+O(x^\th) +O\le(\f{x\log^2 (Q_{\pi}T)}{T}\ri)+\le(\f{\log T}{x}\ri).
\ena
Theorem 5.6 then follows from this and (\ref{perron/NOR+}).
\hfill $\square$

\section{Proof of an almost result}

The following lemma is necessary for Theorem 5.1.

\noindent
{\bf Lemma 5.10.} {\it Let $m\geq 2$ be an integer and let $\pi$ be an irreducible unitary cuspidal
representation of $GL_m({\mathbb A}_{\mathbb Q})$.
Assuming GRH for $L(s,\pi)$, we have
\bea\label{inteXX}
\int^{eX}_X \bigg|\sum_{T<|\g|\leq X^4}\f{x^\rho}{\rho}\bigg|^2
\f{{\rm d} x}{x^2}\ll \f{\log^2 (Q_{\pi}T)}{T},
\eea
for $4\leq T\leq X^4$. }

To prove Lemma 5.10, we need the following lemma of Gallagher
\cite {Gal2}.

\noindent
{\bf Lemma 5.11.} (Gallagher \cite{Gal2}).
{\it  Let
$$
S(u)=\sum_{\nu}c(\nu)e^{2\pi i\nu u}
$$
be absolutely convergent, where the coefficients $c(\nu)\in \mathbb C$, and the
frequencies of $\nu$ run over an arbitrary sequence of real numbers. Then
$$
\int_{-U}^U |S(u)|^2 {\rm d}u
\ll \f{1}{U^2}\int_{-\infty}^{\infty}\bigg|\sum_{x<\nu\leq x+1/U}c(\nu)\bigg|^2{\rm d}x.
$$}

\noindent
{\it Proof of Lemma 5.10.} In the integral in (\ref{inteXX}), we change variables
$x=X e^{2\pi u}$. By GRH, we have $\rho=1/2+i\g$, and therefore
\bea
\int^{eX}_X\bigg|\sum_{T<|\g|\leq X^4}\f{x^\rho}{\rho}\bigg|^2\f{{\rm d}x}{x^2}
&=&2\pi\int^{1/(2\pi)}_0\bigg|\sum_{T<|\g|\leq X^4}
\f{X^{i\g}}{\rho}e^{2\pi i\g u}\bigg|^2 {\rm d}u\nonumber\\
&\ll& \int^1_0\bigg|\sum_{T<|\g|\leq X^4}
\f{X^{i\g}}{\rho}e^{2\pi i\g u}\bigg|^2 {\rm d}u.
\eea
By Gallagher's lemma, the last integral can be estimated as
\bna
&\ll& \int^\infty_{-\infty}
\bigg|\sum_{\substack{T<|\g|\leq X^4\\ t<\g\leq t+1}} \f{X^{i\g}}{\rho}\bigg|^2{\rm d}t
\\
&\ll& \int^\infty_{-\infty}
\bigg\{\sum_{\substack{T<|\g|\leq X^4\\ t<\g\leq t+1}} \f{1}{|\rho|}\bigg\}^2{\rm d}t.
\ena
In the last integral, $t$ should satisfy
either $T-1\leq t\leq X^4$ or $-X^4-1\leq t\leq -T$.
By this and Lemma 4.6,
\bna
\int^\infty_{-\infty}\bigg\{\sum_{\substack{T<|\g|\leq X^4\\ t<\g\leq t+1}}\f{1}{|\rho|}\bigg\}^2{\rm d}t
&\ll&\int^{X^4+1}_{T-1}\bigg\{\sum_{t<\g\leq t+1}\f{1}{|\rho|}\bigg\}^2{\rm d}t
\\
&\ll&\int^{X^4+1}_{T-1}\f{\log^2 (Q_{\pi}t)}{t^2}{\rm d} t
\\\noalign{\vskip 3mm}
&\ll& \f{\log^2 (Q_{\pi}T)}{T}.
\ena
This proves the lemma.
\hfill $\square$

\medskip

Now a proof of Theorem 5.1 is immediate.

\medskip
\noindent
{\it Proof of Theorem 5.1.} Let $2\leq X\leq x\ll X$, and
take $T=X^4$ in the explicit formula (Theorem 5.6). Then
\bea\label{psixpi=}
\psi(x,\pi)=-\sum_{|\g|\leq X^4}\f{x^\rho}{\rho}+O(x^\th \log x)+O\bigg(\f{\log^2(Q_{\pi}x)}{x}\bigg).
\eea
Note that (\ref{psixpi=}) holds unconditionally; while on GRH, the sum
then runs over the nontrivial zeros
$\rho=1/2+i\g$ of $L(s,\pi)$ with $|\g|$ up to $X^4$.

To prove Theorem 5.1, we split the sum over $|\g|$ at $T$, with $2\leq T\leq X^4$
a parameter that will be specified later.

First, we have
\bea\label{psixpi=+}
\sum_{|\g|\leq T}\f{x^\rho}{\rho}
&\ll& x^{1/2}\sum_{|\g|\leq T}\f{1}{|\rho|}
\nonumber
\\
&\ll& x^{1/2}\bigg\{\log Q_\pi+\int_1^T \f{1}{t}{\rm d}N(t,\pi)\bigg\}
\nonumber
\\\noalign{\vskip 2mm}
&\ll& x^{1/2}\log^2 (Q_{\pi}T).
\eea
This inequality asserts that, if $T$ is small, then the contribution to (\ref{psixpi=})
from $|\g|\leq T$ is also small.

However, the contribution to (\ref{psixpi=}) from $T<|\g|\leq X^4$ is not always small; we
shall show that it is usually small. To this end, define
\bna
E(X)=\bigg\{x\in [X,eX]:
\bigg|\sum_{T<|\g|\leq X^4}\f{x^\rho}{\rho}\bigg|
\geq x^{1/2}\log^2(Q_{\pi}\log x)\bigg\}.
\ena
From this and Lemma 5.10, we deduce that
\bna
\log^4(Q_\pi\log X)\int_{E(X)}\f{{\rm d} x}{x}
&\ll&\int_{E(X)}\bigg|\sum_{T<|\g|\leq X^4}\f{x^\rho}{\rho}\bigg|^2\f{{\rm d} x}{x^2}
\\
&\ll&\f{\log^2 (Q_{\pi}T)}{T},
\ena
and hence
\bea\label{psixpi=++}
\int_{E(X)}\f{{\rm d} x}{x} \ll \f{\log^2 (Q_{\pi}T)}{T\log^4(Q_\pi\log X)}.
\eea
Now specify
$$
T=\log X,
$$
and insert (\ref{psixpi=+}) into (\ref{psixpi=}).
Then we see from (\ref{psixpi=++}) that
$$
\psi(x,\pi)\ll x^{1/2}\log^2(Q_{\pi}\log x)
$$
holds on the interval $[X,eX]$ except on the set $E(X)$
of logarithmic measure
$$
\ll \f{1}{T\log^2 (Q_{\pi}T)} \ll \f{1}{T\log^2 T}.
$$

By choosing $X=e^T$ with $T=2,3,\ldots$, the total logarithmic measure
of the exceptional set $E$ will be
$$
\ll \sum_{T=2}^\infty \f{1}{T\log^2 T}<\infty,
$$
and Theorem 5.1 follows. \hfill $\square$

\section{Mean value and Omega estimates for $\psi(x,\pi)$}

\noindent
In this section, we prove Theorems 5.2 and 5.3.

\medskip
\noindent
{\it Proof of Theorem 5.2.} By the explicit formula (\ref{psixpi=}),
\bea\label{146/CH5}
\int^{eX}_X |\psi(x,\pi)|^2\f{{\rm d} x}{x^2}
&\ll& \int^{eX}_X\bigg|\sum_{|\g|\leq X^4}\f{x^\rho}{\rho}\bigg|^2\f{{\rm d} x}{x^2}
+\int^{eX}_X x^{2\th}\log^2 (Q_\pi x)\f{{\rm d} x}{x^2}\nonumber\\
&& +\int^{eX}_X \f{\log^2(Q_{\pi}x)}{x}\f{{\rm d} x}{x^2} .
\eea
The last two integrals are
$\ll \log^2Q_\pi$.
To estimate the first integral on the right-hand side,
we take $T=4$ in Lemma 5.10, getting
\bea
\int^{eX}_X\bigg|\sum_{|\g|\leq X^4}\f{x^\rho}{\rho}\bigg|^2\f{{\rm d} x}{x^2}
&\ll& \log^2 Q_\pi+\int^{eX}_X\bigg|\sum_{4<|\g|\leq X^4}\f{x^\rho}{\rho}\bigg|^2\f{{\rm d} x}{x^2}
\\
&\ll& \log^2 Q_\pi.
\eea
Therefore, (\ref{146/CH5}) becomes
\bea\label{ORG5.3}
\int^{eX}_X |\psi(x,\pi)|^2\f{{\rm d} x}{x^2}\ll \log^2 Q_\pi,
\eea
and consequently,
\bna
\int^{eX}_X |\psi(x,\pi)|^2\f{{\rm d}x}{x}\ll X\log^2 Q_\pi.
\ena
A splitting-up argument then yields
\bna
\int^X_2 |\psi(x,\pi)|\f{{\rm d}x}{x}
&=&\int^{X}_{X/e} |\psi(x,\pi)|^2\f{{\rm d}x}{x}+\int^{X/e}_{X/e^2} |\psi(x,\pi)|^2\f{{\rm d}x}{x}+\cdots
\\
&\ll& \f{X\log^2 Q_\pi}{e} + \f{X\log^2 Q_\pi}{e^2}+\cdots
\\\noalign{\vskip 2mm}
&\ll& X\log^2 Q_\pi.
\ena
This proves Theorem 5.2. \hfill $\square$

\medskip
\noindent
{\it Proof of Theorem 5.3.} To prove the assertion of the
theorem, we assume
\bea\label{5/xalpah}
\psi(x,\pi)\ll x^{\a},
\eea
where $\a$ is some positive absolute constant, and
the $\ll$-constant depends at most on $m$ and $\a$.
We establish the assertion by finding a contradiction whenever
$\a<1/2$.

Applying partial summation to (\ref{5/L'/L}), we get
from the definition of $\psi(x,\pi)$ that,
for $\s>1$,
\bea\label{5/LL=sint}
-\f{L'}{L}(s,\pi)=s\int_{1}^\infty \f{\psi(x,\pi)}{x^{s+1}}{\rm d}x.
\eea
Inserting (\ref{5/xalpah}), we have
$$
\f{\psi(x,\pi)}{x^{s+1}}\ll \f{x^{\a}}{x^{\s+1}}
\ll \f{1}{x^{\s+1-\a}}.
$$
Therefore, for $\s>\a+\ve$ with $\ve>0$ arbitrary, the integral on the
right-hand side of (\ref{5/LL=sint}) is uniformly convergent, and so represents
a regular function in the half-plane $\s>\a$. It follows from
(\ref{5/LL=sint}) that $L(s,\pi)$ cannot have a zero in the half-plane
$\s>\a$. This will lead to a contradiction if $\a<1/2$.
\hfill $\square$

\section{A Linnik-type problem for $\{a_\pi(n)\L(n)\}_{n=1}^\infty$}

\noindent
To each irreducible unitary cuspidal representation
$\pi=\otimes \pi_p$ of $GL_m({\mathbb A}_{\mathbb Q})$, one can attach a global $L$-function
$L(s,\pi)$ as in \S 4.1.
Taking logarithmic differentiation
in (\ref{Lpspip}), one gets (\ref{5/L'/L}) with Dirichlet coefficients
$$
\{a_\pi(n)\L(n)\}_{n=1}^\infty,
$$
where $\L(n)$ is the von Mangoldt function, and
$a_\pi(p^k)$ as in (\ref{apipk/AUTO}), i.e.
\bea\label{apipk/AUCON}
a_\pi(p^k)=
\sum_{j=1}^m
\alpha_\pi(p,j)^k.
\eea

If $\pi$ is self-contragredient, then (\ref{4/for9}) states that
\bna
\{\a_{\pi}(p,j)\}_{j=1}^m
=\{\overline{\a_\pi(p,j)}\}_{j=1}^m,
\ena
and hence, by (\ref{apipk/AUCON}),
\bea\label{realla}
a_{\pi}(p^k)=\overline{a_{\pi}(p^k)},
\eea
which means that $a_{\pi}(n)\L(n)$ are real for all $n\geq 1$.

The purpose of this section is to establish the following
Linnik-type theorem
for the sequence $\{a_\pi(n)\L(n)\}_{n=1}^\infty$.
Actually, this is a corollary to Theorem 4.15 and
L\"u's lemma (Lemma 4.11).

\noindent
{\bf Theorem 5.12.}
{\it Let $m\geq 2$ be an integer and let $\pi$ be an irreducible unitary cuspidal representation for $GL_m(\mathbb A_\mathbb Q)$.
If all $a_\pi(n)\L(n)$ are real, then $\{a_\pi(n)\L(n)\}_{n=1}^\infty$ changes sign at some $n$ satisfying
\bea\label{bound/AUTO-A}
n\ll Q_\pi^{m/2+\ve}.
\eea
The constant implied in (\ref{bound/AUTO-A})
depends only on $m$ and $\ve$. In particular, the result is true
for any self-contragredient representation $\pi$. }

\noindent
{\it Proof.} If we abbreviate $\alpha_\pi(p,j)$ to $\alpha_j$, then
(\ref{apipk/AUCON}) takes the form
\bea\label{apipk/AAA}
a_\pi(p^k)=
\sum_{j=1}^m
\alpha_j^k.
\eea
The key observation is that (\ref{apipk/AAA}) is exactly the same as that defined in (\ref{4/REL/api}),
and Lemma 4.11 is applicable. Thus, as in (\ref{4/FORal++}), we have,
for all $k\geq 1$,
\bea\label{4/FORal++AA}
k\l_{\pi}(p^k)
&=&
a_{\pi}(p)\l_{\pi}(p^{k-1})
+a_{\pi}(p^2)\l_{\pi}(p^{k-2})
+\cdots \nonumber\\\noalign{\vskip 2mm}
&&
+a_{\pi}(p^{k-1})\l_{\pi}(p)
+a_{\pi}(p^k),
\eea
where $\l_\pi(p^k)$ is as in (\ref{lpin=}). By induction on $k$, we show that
if $a_{\pi}(p^k)\geq 0$ for all $k\leq K$, then $\l_{\pi}(p^k)\geq 0$ also for all $k\leq K$.

Now we invoke Theorem 4.15, to deduce that there is an $n$ with
$$
n\ll Q_\pi^{m/2+\ve}
$$
such that $\l_\pi(n)<0$. By (\ref{lpin=}), we see that $\l_\pi(n)$ is multiplicative with respect to $n$,
and therefore there must be a power $p_0^{k_0}$ of a prime $p_0$ with
$$
p_0^{k_0}\ll Q_\pi^{m/2+\ve}
$$
such that $\l_\pi(p_0^{k_0})<0$. Thus, there must be some $k_1\leq k_0$ such that
$a_\pi(p_0^{k_1})<0$. This proves the theorem.
\hfill $\square$

\chapter{Selberg's normal density theorem for automorphic $L$-functions}

\section{Selberg's normal density theorem}

Write, as usual,
$$
\psi(x)=\sum_{n\leq x}\L(n).
$$
It is known that, under the Riemann Hypothesis for the zeta-function,
$$
\psi(x)=x+O(x^{1/2}\log^{2}x).
$$
From this, a result for primes in short intervals
of the form $[x,x+y)$ will follow.

Selberg \cite {Sel} studied
the normal density of primes in short interval. Under the
Riemann Hypothesis for the zeta-function,
i.e. in the case of $m=1$, Selberg \cite {Sel} proved that
\bea\label{6OR1.3}
\int^X_1\{\psi(x+h(x))-\psi(x)-h(x)\}^2 {\rm d}x=o(h(X)^2 X)
\eea
for any increasing functions $h(x) \leq x$ with
$$
\f{h(x)}{\log^2 x}\to \infty.
$$

In Chapter 6, we prove an analogue of this in the case of automorphic $L$-functions attached to
irreducible unitary cuspidal representation
$\pi=\otimes \pi_p$ of $GL_m({\mathbb A}_{\mathbb Q})$.

\section{Selberg's normal density theorem for
automorphic $L$-functions}

To each irreducible unitary cuspidal representation
$\pi=\otimes \pi_p$ of $GL_m({\mathbb A}_{\mathbb Q})$, one can attach a global $L$-function
$L(s,\pi)$ as in \S 4.1.
Let notation be as in the previous chapter.
The prime number theorem for $L(s,\pi)$ concerns the asymptotic behavior of the
counting function
\bna
\psi(x,\pi)=\sum_{n\leq x}\L(n)a_\pi(n),
\ena
and a special case of it
asserts that, if $\pi$ is an irreducible unitary cuspidal
representation of $GL_m({\mathbb A}_{\mathbb Q})$ with $m\geq 2$, then
\bea\label{6OR1.4}
\psi(x,\pi)\ll \sqrt{Q_\pi} \cdot x \cdot \exp\bigg(-\f{c}{2m^4}\sqrt{\log x}\bigg)
\eea
for some absolute positive constant $c$, where the implied constant is absolute.
In Iwaniec and Kowalski
\cite{IwaKow}, Theorem 5.13, a prime number theorem is proved for general
$L$-functions satisfying necessary axioms, from which (\ref{6OR1.4})
follows as a consequence. Under GRH, (\ref{6OR1.4}) can be improved
to
\bea\label{6OR1.5}
\psi(x,\pi)\ll x^{1/2}\log^2 (Q_\pi x).
\eea
It follows from (\ref{6OR1.5}) that, under GRH,
\bea\label{6OR1.6}
\psi(x+h(x),\pi)-\psi(x,\pi)=o(h(x))
\eea
for increasing functions $h(x)\leq x$ satisfying
$$
\f{h(x)}{x^{1/2}\log^2 (Q_\pi x)}\to \infty.
$$
In view of
\bna
\psi(x+h(x),\pi)-\psi(x,\pi)=\sum_{x<n\leq x+h(x)}\L(n)a_\pi (n),
\ena
(\ref{6OR1.6}) describes oscillation of the coefficients $a_\pi(p)$ in the short
intervals $x<p\leq x+h(x)$.

In this Chapter, we will show that (\ref{6OR1.6}) holds on average for even shorter $h(x)$;
see Theorem 6.1 below.

\noindent
{\bf Theorem 6.1.} {\it Let $m\geq 2$ be an integer and let $\pi$ be an irreducible unitary cuspidal
representation of $GL_m({\mathbb A}_{\mathbb Q})$.
Assume GRH for $L(s,\pi)$. We have
\bea\label{OR1.9}
\int^X_1 |\psi(x+h(x),\pi)-\psi(x,\pi)|^2 {\rm d}x=o(h(X)^2 X),
\eea
for any increasing functions $h(x)\leq x$ satisfying
$$
\f{h(x)}{\log^2 (Q_{\pi}x)}\to \infty.
$$
}

Our Theorem 6.1 generalizes Selberg's result to cases when $m\geq 2$.
It also improves an earlier result of the author \cite{Qu2}
that (\ref{OR1.9})
holds for $h(x)\leq x$ satisfying
$$
\f{h(x)}{x^{\th}\log^2 (Q_\pi x)}\to \infty,
$$
where $\th$ is the bound towards the GRC as explained in Lemma 4.8.
The main new idea is a delicate application of Kowalski-Iwaniec's mean value estimate (cf. Lemma 5.9).
We also need an explicit formula established in Chapter 5 in a more precise form.

Unconditionally, Theorem 6.1 would hold for $h(x)=x^{\b}$
with some constant $0<\b<1$. The exact value of $\b$ depends
on two main ingredients: a satisfactory zero-density estimate for
the $L$-function $L(s,\pi)$, and a zero-free region for $L(s,\pi)$
of Littlewood's or Vinogradov's type.

\section{Proof of the theorem}

We prove Theorem 6.1 in this section. The main tools are the explicit formula
in Theorem 5.6 and Lemma 5.10 which is established under GRH.

\noindent
{\it Proof of Theorem 6.1.}
The proof of Theorem 5.6 actually gives an alternative form
of the explicit formula. Let, as in the proof of Theorem 5.6,
$$
-2<a<-1, \qquad b=1+\f{1}{\log x}.
$$
Then the proof of Theorem 5.6 actually gives
\bea\label{6/EXPLICIT}
\psi(x,\pi)
&=&-\sum_{|\g|\leq T}\f{x^\rho}{\rho}
-\sum_{\substack{a<-\l<b\\ |\nu|\leq T}} \f{x^{-\mu}}{-\mu}
+O\bigg\{\sum_{|n-x|\leq x/\sqrt{T}} |\L(n)a_\pi(n)|\bigg\} \nonumber\\
&&
+O\le(\f{x\log^2 (Q_{\pi}x)}{T^{1/2}}\ri) +O\le(\f{\log T}{x}\ri),
\eea
where $\th$ is as in Lemma 4.8, and $\mu$ goes over
the trivial zeros $\mu=\l+i\nu$ of $L(s,\pi)$.

Let $100\leq X\leq x\leq eX$, and take $T=X^4$ in the explicit formula (\ref{6/EXPLICIT}).
Since the length of the interval $(x-x/X^2,x+x/X^2]$ is
$$
\f{2x}{X^2}\leq \f{1}{10},
$$
this interval contains at most one integer; we denote this possible integer by $n_x$.
It follows that
$$
\sum_{|n-x|\leq x/X^2} |\L(n)a_\pi(n)|=|\L(n_x)a_\pi(n_x)|,
$$
and hence (\ref{6/EXPLICIT}) becomes
\bea\label{6/EXPLICIT+}
\psi(x,\pi)
&=&-\sum_{|\g|\leq X^4}\f{x^\rho}{\rho}
-\sum_{\substack{a<-\l<b\\ |\nu|\leq X^4}}\f{x^{-\mu}}{-\mu}
+O\{|\L(n_x)a_\pi(n_x)|+1\} \nonumber\\
&&
+O\le(\f{\log^2 (Q_{\pi}x)}{x}\ri).
\eea
Note that (\ref{6/EXPLICIT+}) holds unconditionally. On GRH, the sum
then runs over the nontrivial zeros
$\rho=1/2+i\g$ of $L(s,\pi)$ with $|\g|$ up to $X^4$. It follows that
\bea\label{6OR4.2}
&& \hskip -10mm
\psi(x+h,\pi)-\psi(x,\pi)\nonumber\\\noalign{\vskip 1mm}
&&\quad
=-\sum_{|\g| \leq X^4}\f{(x+h)^\rho-x^\rho}{\rho}
-\sum_{\substack{a<-\l<b\\ |\nu|\leq X^4}}\f{(x+h)^{-\mu}-x^{-\mu}}{-\mu}\nonumber\\
&&\qquad  +O\{|\L(n_{x+h})a_\pi(n_{x+h})|+|\L(n_x)a_\pi(n_x)|+1\}+O\le(\f{\log^2(Q_{\pi}x)}{x}\ri)
\nonumber\\\noalign{\vskip 1mm}
&&\quad
=: A+B+C+O\le(\f{\log^2(Q_{\pi}x)}{x}\ri),
\eea
say. We start our proof of Theorem 6.1 by estimating the mean value of
(\ref{6OR4.2})
within $X\leq x\leq eX$, while $h$ $(\leq eX)$ is the length of the interval under consideration.
We are interested in how short $h$ can be.

We start from $A$. In $A$, we split the sum over $|\g|$ at $T$, with $4\leq T\leq X^4$
a parameter that will be specified later, and denote
$$
S_1(y,\pi)=\sum_{|\g|\leq T}y^{i\g},
$$
and
$$
S_2(y,\pi)=\sum_{T<|\g|\leq X^4}\f{y^{i\g}}{\rho}.
$$
Then GRH asserts that
\bna
A
&=&
\bigg\{-\sum_{|\g| \leq T}-\sum_{T<|\g| \leq X^4}\bigg\}\f{(x+h)^\rho-x^\rho}{\rho}\\
&=&
-\int^{x+h}_x S_1(y,\pi)\f{dy}{y^{1/2}}+x^{1/2}S_2(x,\pi)-(x+h)^{1/2}S_2(x+h,\pi)
\\\noalign{\vskip 3mm}
&=:&
A_1+A_2+A_3,
\ena
say. Hence,
\bea\label{6OR4.3}
\int_X^{eX} |A|^2{\rm d}x
\ll \int_X^{eX}|A_1|^2{\rm d}x+\int_X^{eX}|A_2|^2{\rm d}x
+\int_X^{eX}|A_3|^2{\rm d}x.
\eea
By Cauchy's inequality,
\bna
|A_1|^2
&\ll& \int^{x+h}_x |S_1(y,\pi)|^2\f{{\rm d}y}{y}\int^{x+h}_x 1^2 {\rm d}y
\\
&=& h\int^{x+h}_x |S_1(y,\pi)|^2\f{{\rm d}y}{y}.
\ena
We note that the upper bound for $h$ is $eX$. Therefore,
the contribution from $|A_1|^2$ is estimated as
\bna
\int_X^{eX}|A_1|^2{\rm d}x
&\ll&
\int_X^{eX}h \int^{x+h}_x |S_1(y,\pi)|^2\f{{\rm d}y}{y}{\rm d}x
\\\noalign{\vskip 2mm}
&\ll&
\int_X^{2eX}h^2|S_1(y,\pi)|^2\f{{\rm d}y}{y}
\\\noalign{\vskip 2mm}
&=&
h^2\int^{2eX}_{X}\bigg|\sum_{|\g|\leq T}y^{i\g}\bigg|^2 \f{{\rm d}y}{y}
\\
&=&
h^2\int^{\log (2e)}_{0}\bigg|\sum_{|\g|\leq T}e^{i\g(u+\log X)}\bigg|^2 {\rm d}u,
\ena
where in the last equality we have changed variables $y=e^{u+\log X}$.
An application of Gallagher's lemma and Lemma 4.6 to the last integral
now leads to
\bna
\int^{\log  (2e)}_{0}\bigg |\sum_{|\g|\leq T}e^{i\g(u+\log X)}\bigg|^2 {\rm d}u
&\ll&
\int^{\infty}_{-\infty}\bigg|\sum_{\substack{|\g|\leq T \\ t<\g\leq t+1}}1\bigg|^2 {\rm d}t
\\
&\ll&
\int^{T}_{0}\bigg\{\sum_{t<\g\leq t+1}1\bigg\}^2 {\rm d}t
\\\noalign{\vskip 3mm}
&\ll& T\log^2 (Q_{\pi}T).
\ena
Thus, we have
\bea\label{6OR4.4}
\int_X^{eX}|A_1|^2{\rm d}x\ll h^2 T\log^2 (Q_{\pi}T).
\eea
The contribution from $|A_2|^2$ can be estimated as
\bea\label{6OR4.5}
\int^{eX}_X|A_2|^2{\rm d}x
&\ll&
X^2\int^{eX}_X|S_2(x,\pi)|^2\f{{\rm d}x}{x}
\nonumber \\
&=& X^2\int^{eX}_X\bigg|\sum_{T<|\g|\leq X^4}
\f{x^{i\g}}{\rho}\bigg|^2\f{{\rm d}x}{x}
\nonumber\\
&\ll& \f{X^2\log^2 (Q_{\pi}T)}{T},
\eea
as a consequence of Lemma 5.10 and GRH. Similarly, after taking $x+h=y$, we have
\bea\label{6OR4.6}
\int^{eX}_X|A_3|^2{\rm d}x
&\ll&
\f{X^2\log^2 (Q_{\pi}T)}{T}.
\eea
We conclude from (\ref{6OR4.3})-(\ref{6OR4.6}) that
\bea\label{6OR4.7}
\int_X^{eX}|A|^2{\rm d}x
\ll
h^2 T\log^2 (Q_{\pi}T)+\f{X^2\log^2 (Q_{\pi}T)}{T}.
\eea

For the mean-value of $|B|^2$, we apply Lemma 4.8, to get
\bea\label{6OR/B+}
\int^{eX}_X|B|^2{\rm d}x
&=& \int^{eX}_X\bigg|\sum_{\substack{a<-\l<b\\ |\nu|\leq X^4}} \f{(x+h)^{-\mu}-x^{-\mu}}{-\mu}\bigg|^2
{\rm d}x
\nonumber\\
&\ll& \int^{eX}_X \bigg\{\sum_{\substack{a<-\l<b\\ |\nu|\leq X^4}} x^{-\l-1}h\bigg\}^2{\rm d}x
\int^{eX}_X (x^{\th-1}h)^2{\rm d}x \nonumber\\
&\ll& X^{2\th-1}h^2\ll h^2.
\eea
It remains to estimate the contribution of $|C|^2$. We have
\bna
\int^{eX}_X|C|^2{\rm d}x
&=&\int^{eX}_X \{|\L(n_{x+h})a_\pi(n_{x+h})|+|\L(n_x)a_\pi(n_x)|+1\}^2{\rm d}x
\\
&\ll& \int^{eX}_X \{|\L(n_{x+h})a_\pi(n_{x+h})|^2
+|\L(n_x)a_\pi(n_x)|^2\} {\rm d}x+X
\\
&\ll &
\sum_{j=[X]}^{[eX]}\int_{j}^{j+1} \{|\L(n_{x+h})a_\pi(n_{x+h})|^2
+|\L(n_x)a_\pi(n_x)|^2\} {\rm d}x +X.
\ena
Since $h(x)$ is increasing and $h(x)\leq x$, we have trivially,
for $j\leq x\leq j+1$, that
$$
j-1\leq n_{x+h(x)}\leq 2(j+2), \quad j-1\leq n_x\leq j+2.
$$
Thus,
\bea\label{6/|C|2}
\int^{eX}_X|C|^2{\rm d}x
&\ll& \sum_{j=[X]-1}^{3[eX]}|\L(j)a_\pi(j)|^2 +X
\nonumber\\
&\ll& X\log^2(Q_\pi X),
\eea
by applying Lemma 5.9.

Finally we apply (\ref{6OR4.7}), (\ref{6OR/B+}), and (\ref{6/|C|2})
to (\ref{6OR4.2}), getting
\bea\label{6/FIN}
\int_{X}^{eX}|\psi(x+h,\pi)-\psi(x,\pi)|^2{\rm d}x
&\ll& h^2 T\log^2 (Q_{\pi}T)+\f{X^2\log^2 (Q_{\pi}T)}{T} \nonumber\\
&& + X\log^2(Q_\pi X)+\f{\log^4(Q_\pi X)}{X}.
\eea
Now we specify the parameter $T$ by taking
$$
h^2 T\log^2 (Q_{\pi}T)=\f{X^2\log^2 (Q_{\pi}T)}{T},
$$
i.e., taking $T=X/h$.
Therefore, the right-hand side of (\ref{6/FIN}) becomes
$$
hX\log^2 (Q_{\pi}X)+X\log^2(Q_\pi X)+\f{\log^4(Q_\pi X)}{X},
$$
which is of order $o(h^2 X)$ as $h\leq eX$ and
$$
\f{h}{\log^2 (Q_{\pi}X)} \to \infty.
$$
Thus for such $h$, we have
\bea\label{6OR4.8}
\int_X^{eX}|\psi(x+h,\pi)-\psi(x,\pi)|^2{\rm d}x
=o(h^2 X).
\eea
In general, let $h=h(x)$ be an increasing function of $x$ satisfying $h(x)\leq x$ and
$$
\f{h(x)}{\log^2 (Q_{\pi}x)} \to \infty.
$$
Then (\ref{6OR4.8}) gives
\bna
\int_{X/e}^{X}|\psi(x+h(x),\pi)-\psi(x,\pi)|^2{\rm d}x
&\ll&
\int_{X/e}^{X}|\psi(x+h(X),\pi)-\psi(x,\pi)|^2{\rm d}x \\
&=& o\le(h(X)^2 \f{X}{e}\ri).
\ena
A splitting-up argument then gives
\bna
\int_1^X |\psi(x+h(x),\pi)-\psi(x,\pi)|^2x
=o(h(X)^2 X),
\ena
and hence our Theorem 6.1.
\hfill $\square$


\end{document}